\documentclass[12pt,twoside]{amsart}
\usepackage{amsmath, amsthm, amscd, amsfonts, amssymb, graphicx}
\usepackage[bookmarksnumbered, plainpages]{hyperref}
\usepackage{marvosym}
\usepackage{geometry}
\makeatletter
\renewcommand\normalsize{%
\@setfontsize\normalsize\@xpt\@xiipt
\abovedisplayskip 5\p@ \@plus1\p@ \@minus6\p@
\abovedisplayshortskip \z@ \@plus3\p@
\belowdisplayshortskip 6\p@ \@plus3\p@ \@minus3\p@
\belowdisplayskip \abovedisplayskip
\let\@listi\@listI}
\makeatother

\newtheorem{thm}{Theorem}[section]

\newtheorem{lem}[thm]{Lemma}

\newtheorem{defn}[thm]{Definition}
\newtheorem{rem}[thm]{\bf{Remark}}

\newtheorem{pf}[thm]{\bf{Proof}}

\numberwithin{equation}{section}

\begin{document}

\title{\bf Left-invariant Ricci collineations associated to canonical connections on three-dimensional Lorentzian Lie groups}
\author{Tao Yu}
\thanks{{\scriptsize
\hskip -0.4 true cm \textit{2010 Mathematics Subject Classification:}
53C40; 53C42.
 }}

\maketitle

\begin{abstract}
In this paper, we classify Left-invariant Ricci collineations associated to canonical connections and Kobayashi-Nomizu connections on three-dimensional Lorentzian Lie groups.

\end{abstract}

\vskip 0.2 true cm


\pagestyle{myheadings}
\markboth{\rightline\scriptsize }
{\leftline{\scriptsize Ricci collineations}}

\bigskip
\bigskip


\section{ Introduction}
\indent The concept of symmetry was firstly proposed by ancient Greek and since the earliest days of natural philosophy, symmetry has furnished insight into the laws of physics and the nature of the cosmos. The two outstanding theoretical achievements of the 20th century, relativity and quantum mechanics, involve notions of symmetry in a fundamental way, being been studied in depth because of their interest from both a mathematical and a physical viewpoint.For example, they simplify Einstein equations and provide a classification of the spacetimes according to the structure of the corresponding Lie algebra,which contains Ricci and curvature collineations ,among others(see \cite{A3},\cite{ A9},\cite{B1},\cite{B2},\cite{B4}).\\
\indent In this article we shall continue to concentrate on Ricci collineations which is, ${\rm L}_{\xi} \widetilde{\rm Ric}^i=0$. In \cite{B5},Michael and Pantelis had done two purposes, the first of which was to present a useful method, which reduced  the computation of the Ricci collineations and the Matter collineations of a given metric to the computation of Killing vectors ,the second of which was  to apply this method and determine all hypersurface orthogonal locally rotationally symmetric  spacetime metrics, which admit proper Matter collineations and proper Ricci collineations. In\cite{C1}, The mathematical concept of left-invariant Ricci collineationsis proposed to denote that  Ricci tensors  has a value of zero with a left-invariant Ricci collineation ~~$\xi$, having the mathematical meaning of a differential homomorphism of preserving Ricci tensors. The purpose of \cite{C1} is to determine all left-invariant Ricci collineations on three-dimensional Lie groups. Since the connection and Ricci tensor are built from the metric tensor, it must inherit its symmetries where any homothetic vector field is a Ricci collineation. In \cite{C2},Wang computed canonical connections and Kobayashi-Nomizu connections and their curvature on three-dimensional Lorentzian Lie groups with a special product structure.\\
\indent In this paper, our motivation is to classify all left-invariant Ricci collineations associated to acanonical connections and Kobayashi-Nomizu connections on three-dimensional Lorentzian Lie group. More precisely,the concept of Ricci collineation was extended to canonical connection and Kobayashi-Nomizu connection on three-deimensional Lorentzian Lie group equiped with a product structure. We classify all left-invariant Ricci collineations on seven different three-dimensional Lie groups.\\
 \indent In section 2, we recall the definition of  the canonical connection and the Kobayashi-Nomizu connection on three-dimensional Lorentzian Lie group with a product structure and their corresponding symmetric tensors and the definition of Ricci collineation with some remarks. In section 3, we classify all left-invariant Ricci collineations associated to canonical connections and Kobayashi-Nomizu connections on three-dimensional unimodular Lorentzian Lie groups. In section 4, we classify all left-invariant Ricci collineations associated to canonical connections and Kobayashi-Nomizu connections on three-dimensional non-unimodular Lorentzian Lie groups.

\vskip 1 true cm
\section{Preliminaries}
\vskip 0.5 true cm
\noindent{\bf 2.1 The symmetric tensors associated to cannonical connection and Kobayashi-Nomizu connection }\\
\indent In this section,we recall the definition of the canonical connection and the Kobayashi-Nomizu connection on three-dimensional Lorentzian Lie group with a product structure and their corresponging symmetric tensors.\\
Let $\nabla$ be the Levi-Civita connection of $G_i$ and $R$ its curvature tensor, taken with the convention
\begin{equation}
R(X,Y)Z=\nabla_X\nabla_YZ-\nabla_Y\nabla_XZ-\nabla_{[X,Y]}Z.\notag
\end{equation}
The Ricci tensor of $(G_i,g)$ is defined by
\begin{equation}{\rm Ric}(X,Y)=-g(R(X,e_1)Y,e_1)-g(R(X,e_2)Y,e_2)+g(R(X,e_3)Y,e_3),\notag
\end{equation}
where $\{e_1,e_2,e_3\}$ is a pseudo-orthonormal basis, with $e_3$ timelike and the Ricci operator Ric is given by
\begin{equation}{\rm Ric}(X,Y)=g({\rm Ric}(X),Y).\notag
\end{equation}
A product structure $J$ on $G_i$ was defined by
\begin{equation}Je_1=e_1,~Je_2=e_2,~Je_3=-e_3,\notag
\end{equation}
then $J^2={\rm id}$ and $g(Je_j,Je_j)=g(e_j,e_j)$. \\
By \cite{A5}, the canonical connection and the Kobayashi-Nomizu connection as follows:
\begin{equation}\nabla^0_XY=\nabla_XY-\frac{1}{2}(\nabla_XJ)JY,\notag
\end{equation}
\begin{equation}\nabla^{{\rm 1}}_XY=\nabla^0_XY-\frac{1}{4}[(\nabla_YJ)JX-(\nabla_{JY}J)X].\notag
\end{equation}
Then their curvature was given by
\begin{equation}
R^0(X,Y)Z=\nabla^0_X\nabla^0_YZ-\nabla^0_Y\nabla^0_XZ-\nabla^0_{[X,Y]}Z,\notag
\end{equation}
\begin{equation}
R^{{1}}(X,Y)Z=\nabla^{{ 1}}_X\nabla^1_YZ-\nabla^1_Y\nabla^1_XZ-\nabla^1_{[X,Y]}Z.\notag
\end{equation}
The Ricci tensors of $(G_i,g)$ associated to the canonical connection and the Kobayashi-Nomizu connection
are defined by
\begin{equation}{\rm Ric}^0(X,Y)=-g(R^0(X,e_1)Y,e_1)-g(R^0(X,e_2)Y,e_2)+g(R^0(X,e_3)Y,e_3),\notag
\end{equation}
\begin{equation}{\rm Ric}^1(X,Y)=-g(R^1(X,e_1)Y,e_1)-g(R^1(X,e_2)Y,e_2)+g(R^1(X,e_3)Y,e_3).\notag
\end{equation}
 The Ricci operators ${\rm Ric}^0$ and ${\rm Ric}^1$ is given by
\begin{equation}{\rm Ric}^0(X,Y)=g({\rm Ric}^0(X),Y),~~{\rm Ric}^1(X,Y)=g({\rm Ric}^1(X),Y).\notag
\end{equation}
Let
\begin{equation}\widetilde{\rm Ric}^0(X,Y)=\frac{{\rm Ric}^0(X,Y)+{\rm Ric}^0(Y,X)}{2},~~\widetilde{\rm Ric}^1(X,Y)=\frac{{\rm Ric}^1(X,Y)+{\rm Ric}^1(Y,X)}{2},\notag
\end{equation}
and
\begin{equation}\widetilde{\rm Ric}^0(X,Y)=g(\widetilde{{\rm Ric}}^0(X),Y),~~\widetilde{\rm Ric}^1(X,Y)=g(\widetilde{{\rm Ric}}^1(X),Y).\notag
\end{equation}
All the analysis so far tells us that $\widetilde{\rm Ric}^0 ~~and ~~\widetilde{\rm Ric}^1$ are symmetric tensors.\\

 Considering some different three-dimensional Lorentzian Lie groups,  we have specific classification about three-dimensional Lorentzian Lie groups in \cite{A2,A4}(see Theorem 2.1 and Theorem 2.2 in \cite{A1}). We shall denote the connected three-dimensional unimodular Lie groups equipped with a left-invariant Lorentzian metric $g$ by $\{G_i\}_{i=1,\cdots,4}$ and having a corresponding Lie algebra $\{\mathfrak{g}\}_{i=1,\cdots,4}$.
\vskip 0.5 true cm
\noindent{\bf 2.2 Three-dimensional unimodular Lorentzian Lie groups }\\
\begin{thm}
Assume~ $(G,g,J)$~was a three-dimensional unimodular Lorentzian Lie group,equipped with a left-invariant Lorentzian metric~$g$~and a product structure~$J$~,and had a Lie algebra ~$\boldsymbol {g}$~ being one of the following:\\
$\boldsymbol {g}_1:$
\begin{equation}
  [e_1,e_2]=\alpha e_1-\beta e_3,~~[e_1,e_3]=-\alpha e_1-\beta e_2,~~[e_2,e_3]=\beta e_1+\alpha e_2+\alpha e_3,~~\alpha\neq0.\notag
\end{equation}
$\boldsymbol {g}_2:$
\begin{equation}
  [e_1,e_2]=\gamma e_2-\beta e_3,~~[e_1,e_3]=-\beta e_2-\gamma e_3,~~[e_2,e_3]=\alpha e_2,~~\gamma\neq0.\notag
\end{equation}
$\boldsymbol {g}_3:$
\begin{equation}
  [e_1,e_2]=-\gamma e_3,~~[e_1,e_3]= e_1-\beta e_2,~~[e_2,e_3]=\alpha e_1,~~\alpha\neq0.\notag
\end{equation}
$\boldsymbol {g}_4:$
\begin{equation}
  [e_1,e_2]=-e_2+(2\eta-\beta)e_3,\eta=1~~or ~~\eta=-1,~~[e_1,e_3]=-\beta e_2+e_3,~~[e_2,e_3]=\alpha e_1,~~.\notag
\end{equation}
\end{thm}
\vskip 0.5 true cm
Next we shall denote the connected three-dimensional non-unimodular Lie groups equipped with a left-invariant Lorentzian metric $g$ by $\{G_i\}_{i=5,6,7}$ and having a corresponding Lie algebra $\{\boldsymbol{g}\}_{i=5,6,7}$.\\
\vskip 0.1 true cm
\noindent{\bf 2.3 Three-dimensional non-unimodular Lorentzian Lie groups }\\
\begin{thm}
Assume~ $(G,g,J)$~was a three-dimensional non-unimodular Lorentzian Lie group,equipped with a left-invariant Lorentzian metric~$g$~and a product structure~$J$~,and had a Lie algebra ~$\boldsymbol {g}$~ being one of the following:\\
$\boldsymbol {g}_5:$
\begin{equation}
  [e_1,e_2]=0,~~[e_1,e_3]=\alpha e_1+\beta e_2,~~[e_2,e_3]=\gamma e_1+\delta e_2,~~\alpha+\delta\neq0~~,\alpha\gamma+\beta\delta=0.\notag
\end{equation}
$\boldsymbol {g}_6:$
\begin{equation}
  [e_1,e_2]=\alpha e_2+\beta e_3,~~[e_1,e_3]=\gamma e_2+\delta e_3,~~[e_2,e_3]=0,~~\alpha+\delta\neq0~~,\alpha\gamma-\beta\delta=0.\notag
\end{equation}
$\boldsymbol {g}_7:$
\begin{equation}
  [e_1,e_2]=-\alpha e_1-\beta e_2-\beta e_3,~~[e_1,e_3]=\alpha e_1+\beta e_2+\beta e_3,~~[e_2,e_3]=\gamma e_1+\delta e_2+\delta e_3,~~\alpha+\delta\neq0~~,\alpha\gamma=0.\notag
\end{equation}
\end{thm}
\begin{defn}\label{2.4}
\vskip 0.5 true cm
A unimodular Lorentzian Lie group of $\{G_i\}_{i=1,\cdots,4}$ admits left-Ricci collineations if and only if it satisfies:\\
\begin{equation}
  {\rm L}_{\xi} \widetilde{\rm Ric}^i=0,
\end{equation}
for~~ $i=1,2$~~,where ~~$\xi~~$ is a element of three-dimensional unimodular Lie group and~~ $\xi=\lambda_1 e_1+\lambda_2 e_2+\lambda_3 e_3$~~,and ~~$\{e_1,e_2,e_3\}$ ~~is a pseudo-orthonormal basis with $e_3$ timelike.
Let
\begin{equation}{\rm L}_{\xi} \widetilde{\rm Ric}^0(X,Y)=\xi(\widetilde{\rm Ric}^0 (X,Y)-\widetilde{\rm Ric}^1([\xi,X],Y)-\widetilde{\rm Ric}^0(X,[\xi,Y]),
\end{equation}
and
\begin{equation}{\rm L}_{\xi} \widetilde{\rm Ric}^1(X,Y)=\xi(\widetilde{\rm Ric}^1(X,Y)-\widetilde{\rm Ric}^1([\xi,X],Y)-\widetilde{\rm Ric}^1(X,[\xi,Y]),
\end{equation}
\end{defn}
\indent According to this definition,we have
\begin{align}\notag
{\rm L}_{e_1} \widetilde{\rm Ric}^i(e_1,e_1)=e_1(\widetilde{\rm Ric}^i (e_1,e_1)-\widetilde{\rm Ric}^i([e_1,e_1],e_1)-\widetilde{\rm Ric}^i(e_1,[e_1,e_1])=0,~~\\\notag
{\rm L}_{e_2} \widetilde{\rm Ric}^i(e_2,e_2)=e_2(\widetilde{\rm Ric}^i (e_2,e_2)-\widetilde{\rm Ric}^i([e_2,e_2],e_2)-\widetilde{\rm Ric}^i(e_2,[e_2,e_2])=0,~~\\\notag
{\rm L}_{e_3} \widetilde{\rm Ric}^i(e_3,e_3)=e_3(\widetilde{\rm Ric}^i (e_3,e_3)-\widetilde{\rm Ric}^i([e_3,e_3],e_3)-\widetilde{\rm Ric}^i(e_3,[e_3,e_3])=0,~~\\\notag
\end{align}
in the following discussion,we will not mention the calculation of these six items ($i=1,2$).owing to $\widetilde{\rm Ric}^0 ~~and~~ \widetilde{\rm Ric}^1$ are symmetric tensors,we have
\begin{align}
{\rm L}_{\xi} \widetilde{\rm Ric}^i (e_1,e_2)={\rm L}_{\xi} \widetilde{\rm Ric}^i (e_2,e_1) ,~~\notag
{\rm L}_{\xi} \widetilde{\rm Ric}^i (e_1,e_3)={\rm L}_{\xi} \widetilde{\rm Ric}^i (e_3,e_1) ,~~\notag
\end{align}
\begin{rem}
 $\mathbb{V}_{\mathbb{R}\mathbb{C}}$ is spanned by Ricci collineations.
\end{rem}
\section{ Left-invariant Ricci collineations associated to canonical connections and Kobayashi-Nomizu connections on three-dimensional unimodular Lorentzian Lie groups}
\vskip 0.5 true cm
\noindent{\bf 3.1 Left-invariant Ricci collineation of $G_1$}\\
\vskip 0.5 true cm
In \cite{A1}, we have for $G_1$, there exists a pseudo-orthonormal basis $\{e_1,e_2,e_3\}$ with $e_3$ timelike such that the Lie
algebra of $G_1$ satisfies
\begin{equation}
[e_1,e_2]=\alpha e_1-\beta e_3,~~[e_1,e_3]=-\alpha e_1-\beta e_2,~~[e_2,e_3]=\beta e_1+\alpha e_2+\alpha e_3,~~\alpha\neq 0.\notag
\end{equation}
By (2.25) in \cite{C2},~~we have
\vskip 0.5 true cm
\begin{lem}
Ricci symmetric tensors of $(G_1,g,J)$  associated to canonical connection are given by
\begin{align}
\widetilde{{\rm Ric}}^0\left(\begin{array}{c}
e_1\\
e_2\\
e_3
\end{array}\right)=\left(\begin{array}{ccc}
-\left(\alpha^2+\frac{\beta^2}{2}\right)&0&-\frac{\alpha\beta}{4}\\
0&-\left(\alpha^2+\frac{\beta^2}{2}\right)&-\frac{\alpha^2}{2}\\
\frac{\alpha\beta}{4}&\frac{\alpha^2}{2}&0
\end{array}\right)\left(\begin{array}{c}
e_1\\
e_2\\
e_3
\end{array}\right).
\end{align}
\end{lem}
By (2.2) and (3.1), we have
\begin{lem}
For $G_1$, the following equalities hold
\begin{align}
&{\rm L}_{e_1} \widetilde{\rm Ric}^0(e_1,e_2)=\alpha(\alpha^2+\frac{3}{4}\beta^2),~~{\rm L}_{e_1} \widetilde{\rm Ric}^0(e_1,e_3)=-\alpha(\alpha^2+\frac{1}{2}\beta^2),~~{\rm L}_{e_1} \widetilde{\rm Ric}^0(e_2,e_2)=\alpha^2 \beta,~~\\\notag
&{\rm L}_{e_1} \widetilde{\rm Ric}^0(e_2,e_3)=-\beta(\frac{5}{4}\alpha^2+\frac{1}{2}\beta^2),~~{\rm L}_{e_1} \widetilde{\rm Ric}^0(e_3,e_3)=\frac{3}{2}\alpha^2 \beta,~~{\rm L}_{e_2} \widetilde{\rm Ric}^0(e_1,e_1)=-2\alpha(\alpha^2+\frac{3}{4}\beta^2)~~\\\notag
&{\rm L}_{e_2} \widetilde{\rm Ric}^0(e_1,e_2)=-\frac{1}{2}\alpha^2 \beta,~~{\rm L}_{e_2} \widetilde{\rm Ric}^0(e_1,e_3)=\beta(\alpha^2+\frac{1}{2}\beta^2)~~,{\rm L}_{e_2} \widetilde{\rm Ric}^0(e_2,e_3)=\frac{1}{2}\alpha(\alpha^2+\beta^2),\\\notag
&{\rm L}_{e_2} \widetilde{\rm Ric}^0(e_3,e_3)=-\alpha(\alpha^2+\frac{1}{2}\beta^2)~~,{\rm L}_{e_3} \widetilde{\rm Ric}^0(e_1,e_1)=2\alpha(\alpha^2+\frac{1}{2}\beta^2),~~{\rm L}_{e_3} \widetilde{\rm Ric}^0(e_1,e_2)=\frac{1}{4}\alpha^2\beta~~\\\notag
&{\rm L}_{e_3} \widetilde{\rm Ric}^0(e_1,e_3)=-\frac{3}{4}\alpha^2 \beta,{\rm L}_{e_3} \widetilde{\rm Ric}^0(e_2,e_2)=-\alpha(\alpha^2+\beta^2)~~,{\rm L}_{e_3} \widetilde{\rm Ric}^0(e_2,e_3)=\frac{1}{2}\alpha(\alpha^2+\frac{1}{2}\beta^2).~~\\\notag
\end{align}
\end{lem}
By (3.2) and definition 2.3, we have
\vskip 0.5 true cm

If the unimodular Lorentzian Lie group $(G_1,g,J)$ admits left-invariant Ricci collineations associated to the canonical connection $\nabla^0$, then
${\rm L}_{\xi} \widetilde{\rm Ric}_{ij}^0=0,for ~~i=1,2,3 ~~and~~ j=1,2,3$,~~ so
\begin{align}
\left\{\begin{array}{l}
-\alpha(\alpha^2+\frac{3}{4}\beta^2)\lambda_2+\alpha(\alpha^2+\frac{1}{2}\beta^2)\lambda_3=0,\\
\\
\alpha(\alpha^2+\frac{3}{4}\beta^2)\lambda_1-\frac{1}{2}\alpha^2 \beta\lambda_2+\frac{1}{4}\alpha^2\beta\lambda_3=0,\\
\\
-\alpha(\alpha^2+\frac{1}{2}\beta^2)\lambda_1+\beta(\alpha^2+\frac{1}{2}\beta^2)\lambda_2-\frac{3}{4}\alpha^2 \beta\lambda_3=0,\\
\\
\alpha^2 \beta\lambda_1-\alpha(\alpha^2+\beta^2)\lambda_3=0,\\
\\
-\beta(\frac{5}{4}\alpha^2+\frac{1}{2}\beta^2)\lambda_1+\frac{1}{2}\alpha(\alpha^2+\beta^2)\lambda_2+\frac{1}{2}\alpha(\alpha^2+\frac{1}{2}\beta^2)\lambda_3=0,\\
\\
\frac{3}{2}\alpha^2 \beta\lambda_1-\alpha(\alpha^2+\frac{1}{2}\beta^2)\lambda_2=0.\\
\end{array}\right.
\end{align}
By solving (3.3), we get
\vskip 0.5 true cm
\begin{thm}
the unimodular Lorentzian Lie group $(G_1,g,J)$ does not admit left-invariant Ricci collineations associated to the canonical connection $\nabla^0$ .
\end{thm}
\begin{pf}
  We analyze each one of these factors by separate. Because $\alpha\neq 0,$~~then we have
 \begin{align}
\left\{\begin{array}{l}
-(\alpha^2+\frac{3}{4}\beta^2)\lambda_2+(\alpha^2+\frac{1}{2}\beta^2)\lambda_3=0,\\
\\
(\alpha^2+\frac{3}{4}\beta^2)\lambda_1-\frac{1}{2}\alpha\beta\lambda_2+\frac{1}{4}\alpha\beta\lambda_2=0,\\
\\
-\alpha(\alpha^2+\frac{1}{2}\beta^2)\lambda_1+\beta(\alpha^2+\frac{1}{2}\beta^2)\lambda_2-\frac{3}{4}\alpha^2\beta\lambda_3=0,\\
\\
\alpha\beta\lambda_1-(\alpha^2+\beta^2)\lambda_3=0,\\
\\
-\beta(\frac{5}{4}\alpha^2+\frac{1}{2}\beta^2)\lambda_1+\frac{1}{2}\alpha(\alpha^2+\beta^2)\lambda_2-\frac{1}{2}\alpha(\alpha^2+\frac{1}{2}\beta^2)\lambda_3=0,\\
\\
\frac{3}{2}\alpha\beta\lambda_1-(\alpha^2+\frac{1}{2}\beta^2)\lambda_2=0.\\
\end{array}\right.
\end{align}
Now if $\beta=0,$~~we have~ $\lambda_1=\lambda_2=\lambda_3=0$.~~If $~~\beta\neq0$,~~naturelly we have $$\lambda_3=\frac{\alpha\beta}{\alpha^2+\beta^2}\lambda_1,\lambda_2=\frac{3\alpha\beta}{2\alpha^2+\beta^2}\lambda_1,$$~~substitute these into Eq.(3.4) to get $$[(\alpha^2+\frac{1}{2}\beta^2)(2\alpha^2+\beta^2)-(\alpha^2+\frac{3}{4}\beta^2)(\alpha^2+\beta^2)]\lambda_1=0,$$~~$$[(4\alpha^2+3\beta^2)(2\alpha^2+\beta^2)(\alpha^2+\beta^2)-\alpha^2\beta^2(4\alpha^2+5\beta^2)]\
\lambda_1=0$$
Solving these, we have $\lambda_1=\lambda_2=\lambda_3=0,\alpha\beta\neq0.$\\
In summary, the unimodular Lorentzian Lie group $(G_1,g,J)$ does not admit left-invariant Ricci collineations associated to the canonical connection $\nabla^0$ .
\end{pf}

\indent By (2.23) in \cite{C2}, we have
\vskip 0.5 true cm
\begin{lem}
Ricci symmetric tensors of $(G_1,g,J)$ ~~associated to Kobayasha-Nomizu connection are given by
\begin{align}
\widetilde{{\rm Ric}}^1\left(\begin{array}{c}
e_1\\
e_2\\
e_3
\end{array}\right)=\left(\begin{array}{ccc}
-\left(\alpha^2+{\beta^2}\right)&\alpha\beta&\frac{\alpha\beta}{2}\\
\alpha\beta&-\left(\alpha^2+{\beta^2}\right)&-\frac{\alpha^2}{2}\\
-\frac{\alpha\beta}{2}&\frac{\alpha^2}{2}&0
\end{array}\right)\left(\begin{array}{c}
e_1\\
e_2\\
e_3
\end{array}\right).
\end{align}
\end{lem}
By (2.3) and (3.5), we have
\begin{lem}
For $G_1$, the following equalities hold
\begin{align}
&{\rm L}_{e_1} \widetilde{\rm Ric}^1(e_1,e_2)=\alpha(\alpha^2+\frac{1}{2}\beta^2),~~{\rm L}_{e_1} \widetilde{\rm Ric}^1(e_1,e_3)=-\alpha^3~~,{\rm L}_{e_1} \widetilde{\rm Ric}^1(e_2,e_2)=-\alpha^2 \beta,\\\notag
&{\rm L}_{e_1} \widetilde{\rm Ric}^1(e_2,e_3)=\beta(\frac{1}{2}\alpha^2-\beta^2)~~,{\rm L}_{e_1} \widetilde{\rm Ric}^1(e_3,e_3)=0,~~{\rm L}_{e_2} \widetilde{\rm Ric}^1(e_1,e_1)=-2\alpha(\alpha^2+\frac{1}{2}\beta^2)~~\\\notag
&{\rm L}_{e_2} \widetilde{\rm Ric}^1(e_1,e_2)=\frac{1}{2}\alpha^2 \beta,~~{\rm L}_{e_2} \widetilde{\rm Ric}^1(e_1,e_3)=\beta^3~~,{\rm L}_{e_2} \widetilde{\rm Ric}^1(e_2,e_3)=\frac{1}{2}\alpha^3,\\\notag
&{\rm L}_{e_2} \widetilde{\rm Ric}^1(e_3,e_3)=\alpha(\beta^2-\alpha^2)~~,{\rm L}_{e_3} \widetilde{\rm Ric}^1(e_1,e_1)=2\alpha^3,~~{\rm L}_{e_3} \widetilde{\rm Ric}^1(e_1,e_2)=-\frac{1}{2}\alpha^2\beta~~\\\notag
&{\rm L}_{e_3} \widetilde{\rm Ric}^1(e_1,e_3)=0,~~{\rm L}_{e_3} \widetilde{\rm Ric}^1(e_2,e_2)=-\alpha^3~~,{\rm L}_{e_3} \widetilde{\rm Ric}^1(e_2,e_3)=\frac{1}{2}\alpha(\alpha^2-\beta^2).~~\\\notag
\end{align}
\end{lem}
By (3.6) and definition 2.3, we have
\vskip 0.3true cm

If the unimodular Lorentzian Lie group $(G_1,g,J)$ admits left-invariant Ricci collineations associated to the Kobayasha-Nomizu connection $\nabla^1$, then
${\rm L}_{\xi} \widetilde{\rm Ric}_{ij}^1=0,for~~i=1,2,3~~ and ~~j=1,2,3$,~~ so
\begin{align}
\left\{\begin{array}{l}
-2\alpha(\alpha^2+\frac{1}{2}\beta^2)\lambda_2+2\alpha^3\lambda_3=0,\\
\\
\alpha(\alpha^2+\frac{1}{2}\beta^2)\lambda_1+\frac{1}{2}\alpha^2\beta\lambda_2-\frac{1}{2}\alpha^2\beta\lambda_3=0,\\
\\
-\alpha^3\lambda_1+\beta^3\lambda_2=0,\\
\\
\alpha^2\beta\lambda_1+\alpha^3\lambda_3=0,\\
\\
\beta(\frac{1}{2}\alpha^2-\beta^2)\lambda_1+\frac{1}{2}\alpha^3\lambda_2+\frac{1}{2}\alpha(\alpha^2-\beta^2)\lambda_3=0,\\
\\
\alpha(\beta^2-\alpha^2)\lambda_2=0.\\
\end{array}\right.
\end{align}
By solving (3.7), we get
\vskip 0.5 true cm
\begin{thm}
the unimodular Lorentzian Lie group ~$(G_1,g,J)$~ does not admit left-invariant Ricci collineations associated to the Kobayashi-Nomizu  connection $\nabla^1$ .
\end{thm}
\begin{pf}
  We analyze each one of these factors by separate. Because $\alpha\neq 0,$~~then we have
 \begin{align}
\left\{\begin{array}{l}
-(\alpha^2+\frac{1}{2}\beta^2)\lambda_2+\alpha^2\lambda_3=0,\\
\\
(\alpha^2+\frac{1}{2}\beta^2)\lambda_1-\frac{1}{2}\alpha\beta\lambda_2-\frac{1}{2}\alpha\beta\lambda_2=0,\\
\\
-\alpha^3\lambda_1+\beta^3\lambda_2=0,\\
\\
\beta\lambda_1+\alpha\lambda_3=0,\\
\\
\beta(\frac{1}{2}\alpha^2-\beta^2)\lambda_1+\frac{1}{2}\alpha^3\lambda_2+\frac{1}{2}\alpha(\alpha^2-\beta^2)\lambda_3=0,\\
\\
(\beta^2-\alpha^2)\lambda_2=0.\\
\end{array}\right.
\end{align}
Now if $\beta^2-\alpha^2=0$,~~we have~ $\alpha=\beta~~or ~~\alpha+\beta=0,$~.If $~\alpha=\beta$~~,naturelly we have $\lambda_1=\lambda_2=-\lambda_3$,substitute these into Eq.(3.8), we get $$-(2\alpha^2+\frac{1}{2}\beta^2)\lambda_2=0,$$
Solving these, we have $\lambda_1=\lambda_2=\lambda_3=0.$\\
If ~~$\beta^2-\alpha^2\neq0,$,~~then we have $\lambda_2=0.$~~Similarly we put this into Eq.(3.8) ,we have $\lambda_2=\lambda_3=0.$\\
In summary,the unimodular Lorentzian Lie group ~$(G_1,g,J)$~ does not admit left-invariant Ricci collineations associated to the Kobayashi-Nomizu connection $\nabla^1$ .
\end{pf}
\vskip 0.5 true cm

\vskip 0.5 true cm
\noindent{\bf 3.2 Left-invariant Ricci collineation of $G_2$}\\
\vskip 0.5 true cm
In \cite{A1}, we have for $G_2$, there exists a pseudo-orthonormal basis $\{e_1,e_2,e_3\}$ with $e_3$ timelike such that the Lie
algebra of $G_2$ satisfies
\begin{equation}
[e_1,e_2]=\gamma e_2-\beta e_3,~~[e_1,e_3]=-\beta e_2-\gamma e_3,~~[e_2,e_3]=\alpha e_1,~~\gamma\neq 0.\notag
\end{equation}
\vskip 0.5 true cm

By (2.44) in \cite{C2}, we have
\vskip 0.5 true cm
\begin{lem}
Ricci symmetric tensors of $(G_2,g,J)$  associated to canonical connection are given by
\begin{align}
\widetilde{{\rm Ric}}^0\left(\begin{array}{c}
e_1\\
e_2\\
e_3
\end{array}\right)=\left(\begin{array}{ccc}
-\left(\gamma^2+\frac{\alpha\beta}{2}\right)&0&0\\
0&-\left(\gamma^2+\frac{\alpha\beta}{2}\right)&\frac{\alpha\gamma}{4}-\frac{\beta\gamma}{2}\\
0&\frac{\beta\gamma}{2}-\frac{\alpha\gamma}{4}&0
\end{array}\right)\left(\begin{array}{c}
e_1\\
e_2\\
e_3
\end{array}\right).
\end{align}
\end{lem}
By (2.2) and (3.9), we have
\begin{lem}
For $G_2$, the following equalities hold
\begin{align}
&{\rm L}_{e_1} \widetilde{\rm Ric}^0(e_1,e_2)=0,~~{\rm L}_{e_1} \widetilde{\rm Ric}^0(e_1,e_3)=0~~,{\rm L}_{e_1} \widetilde{\rm Ric}^0(e_2,e_2)=\gamma(2\gamma^2+\beta^2+\frac{1}{2}\alpha\beta),\\\notag
&{\rm L}_{e_1} \widetilde{\rm Ric}^0(e_2,e_3)=-\beta(\gamma^2+\frac{1}{2}\alpha\beta)~~,{\rm L}_{e_1} \widetilde{\rm Ric}^0(e_3,e_3)=\beta(\beta\gamma-\frac{1}{2}\alpha\gamma),~~{\rm L}_{e_2} \widetilde{\rm Ric}^0(e_1,e_1)=0~~\\\notag
&{\rm L}_{e_2} \widetilde{\rm Ric}^0(e_1,e_2)=-\gamma(\gamma^2+\frac{1}{2}\beta^2+\frac{1}{4}\alpha\beta),~~{\rm L}_{e_2} \widetilde{\rm Ric}^0(e_1,e_3)=\frac{1}{2}\gamma^2(\beta+\frac{3}{2}\alpha)+\frac{1}{2}\alpha^2\beta~~,{\rm L}_{e_2} \widetilde{\rm Ric}^0(e_2,e_3)=0,\\\notag
&{\rm L}_{e_2} \widetilde{\rm Ric}^0(e_3,e_3)=0~~,{\rm L}_{e_3} \widetilde{\rm Ric}^0(e_1,e_1)=0,~~{\rm L}_{e_3} \widetilde{\rm Ric}^0(e_1,e_2)=\frac{1}{2}\gamma^2(\beta-\frac{3}{2}\alpha)+\frac{1}{2}\alpha\beta(\beta-\alpha)~~\\\notag
&{\rm L}_{e_3} \widetilde{\rm Ric}^0(e_1,e_3)=\frac{1}{2}\beta\gamma(\frac{1}{2}\alpha-\beta),~~{\rm L}_{e_3} \widetilde{\rm Ric}^0(e_2,e_2)=0~~,{\rm L}_{e_3} \widetilde{\rm Ric}^0(e_2,e_3)=0.~~\\\notag
\end{align}
\end{lem}

By (3.10) and definition 2.3, we have
\vskip 0.5 true cm

If the unimodular Lorentzian Lie group $(G_2,g,J)$ admits left-invariant Ricci collineations associated to the canonical connection $\nabla^0$, then
${\rm L}_{\xi} \widetilde{\rm Ric}_{ij}^0=0,for ~~i=1,2,3~~ and~~ j=1,2,3$,~~ so
\begin{align}
\left\{\begin{array}{l}
\gamma(\gamma^2+\frac{1}{2}\beta^2+\frac{1}{4}\alpha\beta)\lambda_2-\left[\frac{1}{2}\gamma^2(\beta-\frac{3}{2}\alpha)+\frac{1}{2}\alpha\beta(\beta-\alpha)\right]\lambda_3=0,\\
\\
\left[\frac{1}{2}\gamma^2(\beta+\frac{3}{2}\alpha)+\frac{1}{2}\alpha^2\beta\right]\lambda_2+\frac{1}{2}\beta\gamma(\frac{1}{2}\alpha-\beta)\lambda_3=0,\\
\\
\gamma(2\gamma^2+\beta^2+\frac{1}{2}\alpha\beta)\lambda_1=0,\\
\\
\beta(\gamma^2+\frac{1}{2}\alpha\beta)\lambda_1=0,\\
\\
\beta(\beta\gamma-\frac{1}{2}\alpha\gamma)\lambda_1=0.\\
\end{array}\right.
\end{align}
By solving (3.11), we get
\vskip 0.5 true cm
\begin{thm}
the unimodular Lorentzian Lie group $(G_2,g,J)$ admits left-invariant Ricci collineations associated to the canonical connection $\nabla^0$ if and only if\\
$(1)\alpha=\beta=0,\lambda_1=\lambda_2=0,\mathbb{V}_{\mathbb{R}\mathbb{C}}=\langle e_3\rangle$,\\
\\
$(2)\alpha\neq0,\beta\neq0,\left|\begin{array}{cc}
A&B\\
C&D\\
\end{array}\right|
\neq0,\mathbb{V}_{\mathbb{R}\mathbb{C}}=\langle e_2-\frac{C}{D}e_3\rangle$.
\end{thm}
\begin{pf}
  We analyze each one of these factors by separate. Because $\gamma\neq 0,$~~then we assume that
\begin{align}
\left\{\begin{array}{l}
A=\gamma(2\gamma^2+\beta^2+\frac{1}{2}\alpha\beta),\\
\\
B=\gamma^2(\frac{3}{2}\alpha-\beta)+\alpha\beta(\alpha-\beta),\\
\\
C=\gamma^2(\beta+\frac{3}{2}\alpha)+\alpha^2\beta,\\
\\
D=\beta\gamma(\frac{1}{2}\alpha-\beta).\\
\end{array}\right.
\end{align}
We get$$
\left|\begin{array}{cc}
A&B\\
C&D\\
\end{array}\right|
=4\alpha^4\beta^2-4\alpha^3\beta^3+12\alpha^3\beta\gamma^2-7\alpha^2\beta^2\gamma^2+9\alpha^2\beta^4-4\alpha\beta^3\gamma^2+4\alpha\beta\gamma^4+4\beta^4\gamma^2+4\beta^2\gamma^4
$$
If $\alpha=0$,~~we have$
\left|\begin{array}{cc}
A&B\\
C&D\\
\end{array}\right|
=\beta^2\gamma^2(\beta^2+\gamma^2)
$.\\
Note that if $beta\neq0$~~,$\alpha=0$,~~we have$
\left|\begin{array}{cc}
A&B\\
C&D\\
\end{array}\right|
\neq0
$,then~~$\lambda_1=\lambda_2=\lambda_3=0.$~~\\
Note that  if $beta=0$~~, then$
\left|\begin{array}{cc}
A&B\\
C&D\\
\end{array}\right|
=0$,~~and we have $\lambda_1=\lambda_2=0.$\\
If $\alpha\neq0,\beta=0$,~~we have
 \begin{align}
\left\{\begin{array}{l}
\gamma^3\lambda_2+\frac{3}{4}\alpha\gamma^2\lambda_3=0,\\
\\
\alpha\gamma^2\lambda_2=0,\\
\\
\gamma^3\lambda_1=0.\\
\end{array}\right.
\end{align}
~~Then we have $\lambda_1=\lambda_2=\lambda_3=0.$\\
Notice that when~~$\alpha\neq0,\beta\neq0$,~~we have $\lambda_1=0$,~~whenever $\alpha=2\beta ~~or~~ \alpha \neq2\beta$.
If $\alpha\neq0,\beta\neq0,\left|\begin{array}{cc}
A&B\\
C&D\\
\end{array}\right|
\neq0$,then~~$\lambda_2=\lambda_3=0.$\\
And if $\alpha\neq0,\beta\neq0,\left|\begin{array}{cc}
A&B\\
C&D\\
\end{array}\right|
=0$,then $\lambda_3=-\frac{C}{D}\lambda_2.$\\
In summary,the unimodular Lorentzian Lie group $(G_2,g,J)$ admits left-invariant Ricci collineations associated to the canonical connection $\nabla^0$ if and only if\\
$(1)\alpha=\beta=0,\lambda_1=\lambda_2=0,\mathbb{V}_{\mathbb{R}\mathbb{C}}=\langle e_3\rangle$,\\
$(2)\alpha\neq0,\beta\neq0,\left|\begin{array}{cc}
A&B\\
C&D\\
\end{array}\right|
\neq0,\mathbb{V}_{\mathbb{R}\mathbb{C}}=\langle e_2-\frac{C}{D}e_3\rangle$.
\end{pf}
\indent By (2.54) in \cite{C2}, we have
\vskip 0.5 true cm
\begin{lem}
Ricci symmetric tensors of $(G_2,g,J)$ ~~associated to Kobayasha-Nomizu connection are given by
\begin{align}
\widetilde{{\rm Ric}}^1\left(\begin{array}{c}
e_1\\
e_2\\
e_3
\end{array}\right)=\left(\begin{array}{ccc}
-\left(\beta^2+\gamma^2\right)&0&0\\
0&-\left(\gamma^2+\alpha\beta\right)&\frac{\alpha\gamma}{2}\\
0&-\frac{\alpha\gamma}{2}&0
\end{array}\right)\left(\begin{array}{c}
e_1\\
e_2\\
e_3
\end{array}\right).
\end{align}
\end{lem}
By (2.3) and (3.14),we have
\begin{lem}
For $G_2$, the following equalities hold
\begin{align}
&{\rm L}_{e_1} \widetilde{\rm Ric}^1(e_1,e_2)=0,~~{\rm L}_{e_1} \widetilde{\rm Ric}^1(e_1,e_3)=0~~,{\rm L}_{e_1} \widetilde{\rm Ric}^1(e_2,e_2)=2\gamma(\gamma^2+\frac{1}{2}\alpha\beta),\\\notag
&{\rm L}_{e_1} \widetilde{\rm Ric}^1(e_2,e_3)=-\beta(\gamma^2+\alpha\beta)~~,{\rm L}_{e_1} \widetilde{\rm Ric}^1(e_3,e_3)=-\alpha\beta\gamma,~~{\rm L}_{e_2} \widetilde{\rm Ric}^1(e_1,e_1)=0~~\\\notag
&{\rm L}_{e_2} \widetilde{\rm Ric}^1(e_1,e_2)=-\gamma(\gamma^2+\frac{1}{2}\alpha\beta),~~{\rm L}_{e_2} \widetilde{\rm Ric}^1(e_1,e_3)=\alpha(\beta^2+\frac{1}{2}\gamma^2)~~,{\rm L}_{e_2} \widetilde{\rm Ric}^1(e_2,e_3)=0,\\\notag
&{\rm L}_{e_2} \widetilde{\rm Ric}^1(e_3,e_3)=0~~,{\rm L}_{e_3} \widetilde{\rm Ric}^1(e_1,e_1)=0,~~{\rm L}_{e_3} \widetilde{\rm Ric}^1(e_1,e_2)=\gamma^2(\beta-\frac{1}{2}\alpha)~~\\\notag
&{\rm L}_{e_3} \widetilde{\rm Ric}^1(e_1,e_3)=\frac{1}{2}\alpha\beta\gamma,~~{\rm L}_{e_3} \widetilde{\rm Ric}^1(e_2,e_2)=0~~,{\rm L}_{e_3} \widetilde{\rm Ric}^1(e_2,e_3)=0.~~\\\notag
\end{align}
\end{lem}
By (3.15) and definition 2.3, we have
\vskip 0.1 true cm

If the unimodular Lorentzian Lie group $(G_2,g,J)$ is left-invariant Ricci collineations associated to Kobayashi-Nomizu connection $\nabla^1$, then
${\rm L}_{\xi} \widetilde{\rm Ric}_{ij}^1=0,i=1,2,3~~ and ~~j=1,2,3$, so
\begin{align}
\left\{\begin{array}{l}
\gamma(\gamma^2+\frac{1}{2}\alpha\beta)\lambda_2-\gamma^2(\beta-\frac{1}{2}\alpha)\lambda_3=0,\\
\\
\alpha(\beta^2+\frac{1}{2}\gamma^2)\lambda_2+\frac{1}{2}\alpha\beta\gamma\lambda_3=0,\\
\\
\gamma(\gamma^2+\frac{1}{2}\alpha\beta)\lambda_1=0,\\
\\
\beta(\gamma^2+\alpha\beta)\lambda_1=0,\\
\\
\alpha\beta\gamma\lambda_1=0.\\
\end{array}\right.
\end{align}
By solving (3.16), we get
\vskip 0.5 true cm
\begin{thm}
the unimodular Lorentzian Lie group $(G_2,g,J)$ is left-invariant Ricci collineations associated to the Kobayashi-Nomizu connection $\nabla^1$ if and only if\\
$(1)\alpha=\beta=\lambda_2=\lambda_1=0,\mathbb{V}_{\mathbb{R}\mathbb{C}}=<e_3>,$\\
\\
$(2)\alpha=\lambda_1=0,\beta\neq0,\mathbb{V}_{\mathbb{R}\mathbb{C}}=\langle e_2+\frac{\gamma}{\beta} e_3\rangle,$\\
\\
$(3)\lambda_1=0,\alpha\neq 0,\beta\neq0,\lambda_3=-\frac{2\beta^2+\gamma^2}{\beta\gamma}\lambda_2,\alpha=4\beta,\mathbb{V}_{\mathbb{R}\mathbb{C}}=\langle e_2-\frac{2\beta^2+\gamma^2}{\beta\gamma} e_3\rangle.$\\
\end{thm}
\begin{pf}
  We analyze each one of these factors by separate. Because $\gamma\neq 0,$~~then we have
 \begin{align}
\left\{\begin{array}{l}
(\gamma^2+\frac{1}{2}\alpha\beta)\lambda_2-\gamma(\beta-\frac{1}{2}\alpha)\lambda_3=0,\\
\\
\alpha(\beta^2+\frac{1}{2}\gamma^2)\lambda_2+\frac{1}{2}\alpha\beta\gamma\lambda_3=0,\\
\\
(\gamma^2+\frac{1}{2}\alpha\beta)\lambda_1=0,\\
\\
\beta(\gamma^2+\alpha\beta)\lambda_1=0,\\
\\
\beta(\frac{1}{2}\alpha^2-\beta^2)\lambda_1+\frac{1}{2}\alpha^3\lambda_2+\frac{1}{2}\alpha(\alpha^2-\beta^2)\lambda_3=0,\\
\\
\alpha\beta\lambda_1=0.\\
\end{array}\right.
\end{align}
Notice that compare the third equation with the third fifth equation, $\lambda_1$~~must be equal to 0.\\
If $\alpha=0,$~~we get $\gamma\lambda_2-\beta\lambda_3=0$.~~At the same time,if $\beta=0,$~~then $\lambda_2=0.$~~And if $\beta\neq0,$~~then $\lambda_3=\frac{\gamma}{\beta}\lambda_2$.\\
If $\alpha\neq0,$~~we get
\begin{align}
\left\{\begin{array}{l}
(2\gamma^2+\alpha\beta)\lambda_2+\gamma(\alpha-2\beta)\lambda_3=0,\\
\\
(2\beta^2+\gamma^2)\lambda_2+\beta\gamma\lambda_3=0.\\
\end{array}\right.
\end{align}
Now if $\beta=0,$~~then simply we have $\lambda_2=\lambda_3=0.$\\
we assume that
\begin{align}
\left\{\begin{array}{l}
A=2\gamma^2+\alpha\beta,\\
\\
B=\gamma(\alpha-2\beta),\\
\\
C=2\beta^2+\gamma^2,\\
\\
D=\beta\gamma.\\
\end{array}\right.
\end{align}
We get$$
\left|\begin{array}{cc}
A&B\\
C&D\\
\end{array}\right|
=(\alpha-4\beta)(\beta^2+\gamma^2)
$$
Now if $(\alpha-4\beta)(\beta^2+\gamma^2)\neq0$,equal to $\alpha\neq4\beta$,~~we have~ $\lambda_2=\lambda_3=0$.~~If $(\alpha-4\beta)(\beta^2+\gamma^2)=0$, equal to $\alpha=4\beta$,~~we have~ $\lambda_3=-\frac{2\beta^2+\gamma^2}{\beta\gamma}\lambda_2$.\\
In summary, the unimodular Lorentzian Lie group $(G_2,g,J)$ is left-invariant Ricci collineations associated to the Kobayashi-Nomizu connection $\nabla^1$ if and only if\\
$(1)\alpha=\beta=\lambda_2=\lambda_1=0,\mathbb{V}_{\mathbb{R}\mathbb{C}}=<e_3>,$\\
\\
$(2)\alpha=\lambda_1=0,\beta\neq0,\mathbb{V}_{\mathbb{R}\mathbb{C}}=\langle e_2+\frac{\gamma}{\beta} e_3\rangle,$\\
\\
$(3)\lambda_1=0,\alpha\neq 0,\beta\neq0,\lambda_3=-\frac{2\beta^2+\gamma^2}{\beta\gamma}\lambda_2,\alpha=4\beta,\mathbb{V}_{\mathbb{R}\mathbb{C}}=\langle e_2-\frac{2\beta^2+\gamma^2}{\beta\gamma} e_3\rangle.$\\
\end{pf}
\vskip 1 true cm

\vskip 0.5 true cm
\noindent{\bf 3.3 Left-invariant Ricci collineation of $G_3$}\\
\vskip 0.5 true cm
By \cite{A1}, we have for $G_3$, there exists a pseudo-orthonormal basis $\{e_1,e_2,e_3\}$ with $e_3$ timelike such that the Lie
algebra of $G_3$ satisfies
\begin{equation}
[e_1,e_2]=\gamma e_3,~~[e_1,e_3]=-\beta e_2,~~[e_2,e_3]=\alpha e_1,~~.\notag
\end{equation}
\vskip 0.2 true cm
By (2.64) in \cite{C2}, we have
\vskip 0.5 true cm
\begin{lem}
Ricci symmetric tensors of $(G_3,g,J)$  associated to canonical connection are given by
\begin{align}
\widetilde{{\rm Ric}}^0\left(\begin{array}{c}
e_1\\
e_2\\
e_3
\end{array}\right)=\left(\begin{array}{ccc}
-\gamma a_3&0&0\\
0&-\gamma a_3&0\\
0&0&0
\end{array}\right)\left(\begin{array}{c}
e_1\\
e_2\\
e_3
\end{array}\right).
\end{align}
\end{lem}
where
\begin{equation}
  a_1=\frac{1}{2}(\alpha-\beta-\gamma),a_2=\frac{1}{2}(\alpha-\beta+\gamma),a_3=\frac{1}{2}(\alpha+\beta-\gamma)\notag
\end{equation}
By (2.2) and (3.20), we have
\begin{lem}
For $G_2$, the following equalities hold
\vskip 0.1 true cm
\begin{align}
&{\rm L}_{e_1} \widetilde{\rm Ric}^0(e_1,e_2)=0,~~{\rm L}_{e_1} \widetilde{\rm Ric}^0(e_1,e_3)=0~~,{\rm L}_{e_1} \widetilde{\rm Ric}^0(e_2,e_2)=0,\\\notag
&{\rm L}_{e_1} \widetilde{\rm Ric}^0(e_2,e_3)=-\beta\gamma a_3~~,{\rm L}_{e_1} \widetilde{\rm Ric}^0(e_3,e_3)=0,~~{\rm L}_{e_2} \widetilde{\rm Ric}^0(e_1,e_1)=0~~\\\notag
&{\rm L}_{e_2} \widetilde{\rm Ric}^0(e_1,e_2)=0,~~{\rm L}_{e_2} \widetilde{\rm Ric}^0(e_1,e_3)=\alpha\gamma a_3~~,{\rm L}_{e_2} \widetilde{\rm Ric}^0(e_2,e_3)=0,\\\notag
&{\rm L}_{e_2} \widetilde{\rm Ric}^0(e_3,e_3)=0~~,{\rm L}_{e_3} \widetilde{\rm Ric}^0(e_1,e_1)=0,~~{\rm L}_{e_3} \widetilde{\rm Ric}^0(e_1,e_2)=(\beta-\alpha)\gamma a_3~~\\\notag
&{\rm L}_{e_3} \widetilde{\rm Ric}^0(e_1,e_3)=0,~~{\rm L}_{e_3} \widetilde{\rm Ric}^0(e_2,e_2)=0~~,{\rm L}_{e_3} \widetilde{\rm Ric}^0(e_2,e_3)=0.~~\\\notag
\end{align}
\end{lem}
By (3.21) and definition 2.3, we have
\vskip 0.5 true cm

If the unimodular Lorentzian Lie group $(G_3,g,J)$ admits left-invariant Ricci collineations associated to the canonical connection $\nabla^0$, then
${\rm L}_{\xi} \widetilde{\rm Ric}_{ij}^0=0,i=1,2,3~~ and~~ j=1,2,3$, so
\begin{align}
\left\{\begin{array}{l}
(\beta-\alpha)\gamma a_3\lambda_3=0,\\
\\
\alpha\gamma a_3\lambda_2=0,\\
\\
\beta\gamma a_3\lambda_1=0,\\
\end{array}\right.
\end{align}
By solving (3.22), we get
\vskip 0.5 true cm
\begin{thm}
the unimodular Lorentzian Lie group $(G_3,g,J)$ admits left-invariant Ricci collineations associated to the canonical connection $\nabla^0$ if and only if\\
$(1)\gamma=0,\mathbb{V}_{\mathbb{R}\mathbb{C}}=\langle e_1,e_2,e_3\rangle,$\\
\\
$(2)\gamma\neq0,\alpha=\beta=0,\mathbb{V}_{\mathbb{R}\mathbb{C}}=\langle e_1,e_2,e_3\rangle,$\\
\\
$(3)\alpha=0,\gamma\neq0,\beta\neq0,\beta=\gamma,\mathbb{V}_{\mathbb{R}\mathbb{C}}=\langle e_1,e_2,e_3\rangle,$\\
\\
$(4)\alpha=0,\gamma\neq0,\beta\neq0,\beta\neq\gamma,\mathbb{V}_{\mathbb{R}\mathbb{C}}=\langle e_2\rangle,$\\
\\
$(5)\gamma\neq0,\alpha\neq0,\gamma=\alpha+\beta,\mathbb{V}_{\mathbb{R}\mathbb{C}}=\langle e_1,e_2,e_3\rangle,$\\
\\
$(6)\gamma\neq0,\alpha\neq0,\gamma\neq\alpha+\beta,\alpha=\beta,\mathbb{V}_{\mathbb{R}\mathbb{C}}=\langle e_3\rangle.$\\
\end{thm}
\begin{pf}
  We analyze each one of these factors by separate. Because $a_3=\frac{1}{2}(\alpha+\beta-\gamma),$~~then we have
 \begin{align}
\left\{\begin{array}{l}
\gamma(\beta-\alpha)(\alpha+\beta-\gamma)\lambda_3=0,\\
\\
\alpha\gamma\lambda_2(\alpha+\beta-\gamma)=0,\\
\\
\beta\gamma\lambda_1(\alpha+\beta-\gamma)=0,\\
\end{array}\right.
\end{align}
If $\gamma=0,$~~we get $\lambda_1\in R,\lambda_2\in R,\lambda_3\in R$.~~\\
If $\gamma\neq0,$~~we get
\begin{align}
\left\{\begin{array}{l}
(\beta-\alpha)\lambda_3(\alpha+\beta-\gamma)=0,\\
\\
\alpha\lambda_2(\alpha+\beta-\gamma)=0,\\
\\
\beta\lambda_1(\alpha+\beta-\gamma)=0.\\
\end{array}\right.
\end{align}
Now if $\alpha=0,\beta=0$~~,then simply we have $\lambda_1\in R,\lambda_2\in R,\lambda_3\in R.$\\
If~~$\alpha=0,\beta\neq0$,~~we have
\begin{align}
\left\{\begin{array}{l}
(\beta-\gamma)\lambda_3=0,\\
\\
(\beta-\gamma)\lambda_1=0.\\
\end{array}\right.
\end{align}
then if $\beta-\gamma=0$,~~simply having $\lambda_1\in R,\lambda_2\in R,\lambda_3\in R$~~and if $\beta-\gamma\neq0$,~~simply having $\lambda_2=\lambda_3=0.$\\
Now if $\alpha\neq0,$,~~we have
\begin{align}
\left\{\begin{array}{l}
(\beta-\alpha)(\alpha+\beta-\gamma)\lambda_3=0,\\
\\
(\alpha+\beta-\gamma)\lambda_2=0,\\
\\
\beta(\alpha+\beta-\gamma)\lambda_1=0.\\
\end{array}\right.
\end{align}
In this case, if $\gamma=\alpha+\beta,$~~we have $\lambda_1\in R,\lambda_2\in R,\lambda_3\in R$~~ and if $\gamma\neq\alpha+\beta,$~~we have $\lambda_2=0$.~~Furthermore if $\beta\neq0,\alpha\neq\beta,$~~we get $\lambda_1=\lambda_3=0,$~~and if $\beta\neq0,\alpha=\beta,$~~we get $\lambda_1=0 $~~where if $\beta=0,$~~then we have $\lambda_3=0.$\\
In summary,the unimodular Lorentzian Lie group $(G_3,g,J)$ admits left-invariant Ricci collineations associated to the canonical connection $\nabla^0$ if and only if\\
$(1)\gamma=0,\mathbb{V}_{\mathbb{R}\mathbb{C}}=\langle e_1,e_2,e_3\rangle,$\\
\\
$(2)\gamma\neq0,\alpha=\beta=0,\mathbb{V}_{\mathbb{R}\mathbb{C}}=\langle e_1,e_2,e_3\rangle,$\\
\\
$(3)\alpha=0,\gamma\neq0,\beta\neq0,\beta=\gamma,\mathbb{V}_{\mathbb{R}\mathbb{C}}=\langle e_1,e_2,e_3\rangle,$\\
\\
$(4)\alpha=0,\gamma\neq0,\beta\neq0,\beta\neq\gamma,\mathbb{V}_{\mathbb{R}\mathbb{C}}=\langle e_2\rangle,$\\
\\
$(5)\gamma\neq0,\alpha\neq0,\gamma=\alpha+\beta,\mathbb{V}_{\mathbb{R}\mathbb{C}}=\langle e_1,e_2,e_3\rangle,$\\
\\
$(6)\gamma\neq0,\alpha\neq0,\gamma\neq\alpha+\beta,\alpha=\beta,\mathbb{V}_{\mathbb{R}\mathbb{C}}=\langle e_3\rangle.$\\
\end{pf}
\indent By (2.69) in \cite{C2}, we have
\vskip 0.1 true cm
\begin{lem}
Ricci symmetric tensors of $(G_3,g,J)$ ~~associated to Kobayashi-Nomizu connection are given by
\begin{align}
\widetilde{{\rm Ric}}^1\left(\begin{array}{c}
e_1\\
e_2\\
e_3
\end{array}\right)=\left(\begin{array}{ccc}
\gamma(a_1-a_3)&0&0\\
0&-\gamma(a_2+a_3)&0\\
0&0&0
\end{array}\right)\left(\begin{array}{c}
e_1\\
e_2\\
e_3
\end{array}\right).
\end{align}
\end{lem}
By (2.3) and (3.27), we have
\begin{lem}
For $G_3$, the following equalities hold
\begin{align}
&{\rm L}_{e_1} \widetilde{\rm Ric}^1(e_1,e_2)=0,~~{\rm L}_{e_1} \widetilde{\rm Ric}^1(e_1,e_3)=0~~,{\rm L}_{e_1} \widetilde{\rm Ric}^1(e_2,e_2)=0,\\\notag
&{\rm L}_{e_1} \widetilde{\rm Ric}^1(e_2,e_3)=-\beta\gamma(a_2+a_3)~~,{\rm L}_{e_1} \widetilde{\rm Ric}^1(e_3,e_3)=0,~~{\rm L}_{e_2} \widetilde{\rm Ric}^1(e_1,e_1)=0~~\\\notag
&{\rm L}_{e_2} \widetilde{\rm Ric}^1(e_1,e_2)=0,~~{\rm L}_{e_2} \widetilde{\rm Ric}^1(e_1,e_3)=\alpha\gamma(a_3-a_1)~~,{\rm L}_{e_2} \widetilde{\rm Ric}^1(e_2,e_3)=0,\\\notag
&{\rm L}_{e_2} \widetilde{\rm Ric}^1(e_3,e_3)=0~~,{\rm L}_{e_3} \widetilde{\rm Ric}^1(e_1,e_1)=0,~~{\rm L}_{e_3} \widetilde{\rm Ric}^1(e_1,e_2)=\beta\gamma(a_2+a_3)+\alpha\gamma(a_1-a_3)~~\\\notag
&{\rm L}_{e_3} \widetilde{\rm Ric}^1(e_1,e_3)=0,~~{\rm L}_{e_3} \widetilde{\rm Ric}^1(e_2,e_2)=0~~,{\rm L}_{e_3} \widetilde{\rm Ric}^1(e_2,e_3)=0.~~\\\notag
\end{align}
\end{lem}
By (3.28) and definition 2.3, we get
\vskip 0.3 true cm
If the unimodular Lorentzian Lie group $(G_3,g,J)$ admits left-invariant Ricci collineations associated to the Kobayashi-Nomizu connection $\nabla^1$, then
${\rm L}_{\xi} \widetilde{\rm Ric}_{ij}^1=0,i=1,2,3~~ and ~~j=1,2,3$, so
\begin{align}
\left\{\begin{array}{l}
[\beta\gamma(a_2+a_3)+\alpha\gamma(a_1-a_3)]\lambda_3=0,\\
\\
\alpha\gamma(a_3-a_1)\lambda_2=0,\\
\\
\beta\gamma(a_2+a_3)\lambda_1=0.\\
\end{array}\right.
\end{align}
By solving (3.29), we get
\vskip 0.5 true cm
\begin{thm}
the unimodular Lorentzian Lie group $(G_3,g,J)$ is left-invariant Ricci collineations associated to the Kobayashi-Nomizu connection $\nabla^1$ if and only if ( easy to prove )\\
$(1)\alpha\beta\gamma=0,\mathbb{V}_{\mathbb{R}\mathbb{C}}=\langle e_1,e_2, e_3\rangle,$~~\\
\\
$(2)\alpha\beta\gamma\neq0,\mathbb{V}_{\mathbb{R}\mathbb{C}}=\langle e_3\rangle.$\\
\end{thm}
\vskip 0.5 true cm
\noindent{\bf 3.4 Left-invariant Ricci collineation of $G_4$}\\
\vskip 0.5 true cm
By  \cite{A1}, we have for $G_4$, there exists a pseudo-orthonormal basis $\{e_1,e_2,e_3\}$ with $e_3$ timelike such that the Lie
algebra of $G_4$ satisfies
\begin{equation}
[e_1,e_2]=-e_2+(2\eta-\beta)e_3,\eta=1~~or~~\eta=-1,~~[e_1,e_3]=-\beta e_2+e_3,~~[e_2,e_3]=\alpha e_1~~.\notag
\end{equation}
\vskip 0.5 true cm

By (2.81) in \cite{C2}, we have
\vskip 0.5 true cm
\begin{lem}
Ricci symmetric tensors of $(G_4,g,J)$  associated to canonical connection are given by
\begin{align}
\widetilde{{\rm Ric}}^0\left(\begin{array}{c}
e_1\\
e_2\\
e_3
\end{array}\right)=\left(\begin{array}{ccc}
\left[(2\eta-\beta)b_3-1\right]&0&0\\
0&\left[(2\eta-\beta)b_3-1\right]&\frac{\beta-b_3}{2}\\
0&\frac{b_3-\beta}{2}&0
\end{array}\right)\left(\begin{array}{c}
e_1\\
e_2\\
e_3
\end{array}\right).
\end{align}
\end{lem}
where
\begin{equation}
 b_1=\frac{\alpha}{2}+\eta-\beta,b_2=\frac{\alpha}{2}-\eta,b_3=\frac{\alpha}{2}+\eta\notag
\end{equation}
By (3.30) and (2.3)
\begin{lem}
For $G_4$, the following equalities hold
\vskip 0.1 true cm
\begin{align}
&{\rm L}_{e_1} \widetilde{\rm Ric}^0(e_1,e_2)=0,~~{\rm L}_{e_1} \widetilde{\rm Ric}^0(e_1,e_3)=0~~,{\rm L}_{e_1} \widetilde{\rm Ric}^0(e_2,e_2)=(2\eta-\beta)(b_3+\beta)-2,\\\notag
&{\rm L}_{e_1} \widetilde{\rm Ric}^0(e_2,e_3)=\beta[(2\eta-\beta)b_3-1]~~,{\rm L}_{e_1} \widetilde{\rm Ric}^0(e_3,e_3)=\beta(b_3-\beta),~~{\rm L}_{e_2} \widetilde{\rm Ric}^0(e_1,e_1)=0~~\\\notag
&{\rm L}_{e_2} \widetilde{\rm Ric}^0(e_1,e_2)=(b_3+\beta)(\frac{1}{2}\beta-\eta)+1,~~{\rm L}_{e_2} \widetilde{\rm Ric}^0(e_1,e_3)=\alpha[(\beta-2\eta)b_3+1]+\frac{1}{2}(\beta-b_3)~~,{\rm L}_{e_2} \widetilde{\rm Ric}^0(e_2,e_3)=0,\\\notag
&{\rm L}_{e_2} \widetilde{\rm Ric}^0(e_3,e_3)=0~~,{\rm L}_{e_3} \widetilde{\rm Ric}^0(e_1,e_1)=0,~~{\rm L}_{e_3} \widetilde{\rm Ric}^0(e_1,e_2)=(\alpha-\beta)[(2\eta-\beta)b_3-1]+\frac{1}{2}(b_3-\beta)~~,\\\notag
&{\rm L}_{e_3} \widetilde{\rm Ric}^0(e_1,e_3)=\frac{1}{2}\beta(\beta-b_3),~~{\rm L}_{e_3} \widetilde{\rm Ric}^0(e_2,e_2)=0~~,{\rm L}_{e_3} \widetilde{\rm Ric}^0(e_2,e_3)=0.~~\\\notag
\end{align}
\end{lem}
By (3.31) and definition 2.3, we have
\vskip 0.5 true cm

If the unimodular Lorentzian Lie group $(G_4,g,J)$ admits left-invariant Ricci collineations associated to the canonical connection $\nabla^0$, then
${\rm L}_{\xi} \widetilde{\rm Ric}_{ij}^0=0,for~~i=1,2,3~~ and~~ j=1,2,3$,~~ so
\begin{align}
\left\{\begin{array}{l}
\left[(b_3+\beta)(\frac{1}{2}\beta-\eta)+1\right]\lambda_2+\{(\alpha-\beta)\left[(2\eta-\beta)b_3-1\right]+\frac{1}{2}(b_3-\beta)\}\lambda_3=0,\\
\\
\{\alpha[(\beta-2\eta)b_3+1]+\frac{1}{2}(\beta-b_3)\}\lambda_2+\frac{1}{2}\beta(\beta-b_3)\lambda_3=0,\\
\\
\left[(2\eta-\beta)(b_3+\beta)-2\right]\lambda_1=0,\\
\\
\beta\left[(2\eta-\beta)b_3-1\right]\lambda_1=0,\\
\\
\beta(b_3-\beta)\lambda_1=0.\\
\end{array}\right.
\end{align}
By solving (3.32), we get
\vskip 0.5 true cm
\begin{thm}
the unimodular Lorentzian Lie group ~$(G_4,g,J)$~ admits left-invariant Ricci collineations associated to the canonical connection $\nabla^0$ if and only if\\
$(1)\beta=0,\eta=1,\alpha=0,\lambda_2=\lambda_3=0,\mathbb{V}_{\mathbb{R}\mathbb{C}}=\langle e_1\rangle,$\\
\\
$(2)\beta=0,\eta=1,(\alpha+2)(\alpha+\frac{1}{4})=0,\lambda_1=\lambda_2=0,\mathbb{V}_{\mathbb{R}\mathbb{C}}=\langle e_3\rangle,$\\
\\
$(3)\beta=0,\eta=-1,\alpha=0,\lambda_2=\lambda_3=0,\mathbb{V}_{\mathbb{R}\mathbb{C}}=\langle e_1\rangle,$\\
\\
$(4)\beta=0,\eta=-1,(\alpha-2)(\alpha-\frac{1}{4})=0,\lambda_1=\lambda_2=0,\mathbb{V}_{\mathbb{R}\mathbb{C}}=\langle e_3\rangle,$\\
\\
$(5)\beta\neq0,\eta=1,\alpha=0,\beta=1,\mathbb{V}_{\mathbb{R}\mathbb{C}}=\langle e_1,e_2,e_3\rangle,$\\
\\
$(6)\beta\neq0,\eta=1,\alpha=2,\beta=2,\mathbb{V}_{\mathbb{R}\mathbb{C}}=\langle e_3\rangle,$\\
\\
$(7)\beta\neq0,\eta=1,\alpha+2\neq2\beta,
\left|\begin{array}{cc}
A&B\\
C&D\\
\end{array}\right|
=0,
\mathbb{V}_{\mathbb{R}\mathbb{C}}=\langle e_2-\frac{C}{D} e_3\rangle,$\\
\\
$(8)\beta\neq0,\eta=-1,\alpha=0,\beta=-1,\mathbb{V}_{\mathbb{R}\mathbb{C}}=\langle e_1,e_2,e_3\rangle,$\\
\\
$(9)\beta\neq0,\eta=-1,\alpha=-2,\beta=-2,\mathbb{V}_{\mathbb{R}\mathbb{C}}=\langle e_3\rangle,$\\
\\
$(10)\beta\neq0,\eta=-1,\alpha-2\neq2\beta,
\left|\begin{array}{cc}
A&B\\
C&D\\
\end{array}\right|
=0,
\mathbb{V}_{\mathbb{R}\mathbb{C}}=\langle e_2-\frac{C}{D}e_3\rangle.$\\
\\
\end{thm}
\begin{pf}
  We analyze each one of these factors by separate. Because $b_3=\frac{1}{2}\alpha+\beta,$~~if $\beta=0,$~~then we have
 \begin{align}
\left\{\begin{array}{l}
(2\eta b_3-2)\lambda_1=0,\\
\\
(-b_3\eta+1)\lambda_2+(2\alpha\eta b_3+\frac{1}{2}b_3)\lambda_3=0,\\
\\
(-2\alpha\eta b_3+\alpha-\frac{1}{2}b_3)\lambda_2=0.\\
\end{array}\right.
\end{align}
Now if $\eta=1,$~~we have
\begin{align}
\left\{\begin{array}{l}
\alpha\lambda_1=0,\\
\\
-\frac{\alpha}{2}\lambda_2+(\alpha^2+\frac{9}{4}\alpha+\frac{1}{2})\lambda_3=0,\\
\\
(\alpha^2+\frac{5}{4}\alpha+\frac{1}{2})\lambda_2=0.\\
\end{array}\right.
\end{align}
Furthermore if $\alpha=0,$~~then we get $\lambda_2=\lambda_3=0.$~~If $\alpha\neq0$~~ and $(\alpha+2)(\alpha+\frac{1}{4})\neq0,$~~we get~~$\lambda_1=\lambda_2=\lambda_3=0.$~~If $\alpha\neq0$~~ and $(\alpha+2)(\alpha+\frac{1}{4})=0,$~~we get~~$\lambda_1=\lambda_2=0.$~~\\
Although we assume $\eta=-1,$~~we still have
\begin{align}
\left\{\begin{array}{l}
\alpha\lambda_1=0,\\
\\
\frac{\alpha}{2}\lambda_2+(-\alpha^2+\frac{9}{4}\alpha-\frac{1}{2})\lambda_3=0,\\
\\
(\alpha^2-\frac{5}{4}\alpha+\frac{1}{2})\lambda_2=0.\\
\end{array}\right.
\end{align}
Furthermore if $\alpha=0,$~~then we get $\lambda_2=\lambda_3=0.$~~If ~~$\alpha\neq0$~~ and $(\alpha-2)(\alpha-\frac{1}{4})\neq0,$~~we get~~$\lambda_1=\lambda_2=\lambda_3=0.$~~If $\alpha\neq0$~~ and $(\alpha-2)(\alpha-\frac{1}{4})=0,$~~we~~get~~$\lambda_1=\lambda_2=0.$~~\\
If~$\beta\neq0,$~~then we have
 \begin{align}
\left\{\begin{array}{l}
\left[(b_3+\beta)(\frac{1}{2}\beta-\eta)+1\right]\lambda_2+\{(\alpha-\beta)\left[(2\eta-\beta)b_3-1\right]+\frac{1}{2}(b_3-\beta)\}\lambda_3=0,\\
\\
\{\alpha[(\beta-2\eta)b_3+1]+\frac{1}{2}(\beta-b_3)\}\lambda_2+\frac{1}{2}\beta(\beta-b_3)\lambda_3=0,\\
\\
\left[(2\eta-\beta)(b_3+\beta)-2\right]\lambda_1=0,\\
\\
\left[(2\eta-\beta)b_3-1\right]\lambda_1=0,\\
\\
(b_3-\beta)\lambda_1=0.\\
\end{array}\right.
\end{align}
Now if $\eta=1,$~~we have
\begin{align}
\left\{\begin{array}{l}
\left[(\frac{\alpha}{2}+1+\beta)(\frac{1}{2}\beta-1)+1\right]\lambda_2+\{(\alpha-\beta)\left[(2-\beta)(\frac{\alpha}{2}+1)-1\right]+\frac{1}{2}(\frac{\alpha}{2}+1-\beta)\}\lambda_3=0,\\
\\
\{\alpha[(\beta-2)(\frac{\alpha}{2}+1)+1]+\frac{1}{2}(\beta-\frac{\alpha}{2}-1)\}\lambda_2+\frac{1}{2}\beta(\beta-\frac{\alpha}{2}-1)\lambda_3=0,\\
\\
\left[(2-\beta)(\frac{\alpha}{2}+1+\beta)-2\right]\lambda_1=0,\\
\\
\left[(2-\beta)(\frac{\alpha}{2}+1)-1\right]\lambda_1=0,\\
\\
(\frac{\alpha}{2}+1-\beta)\lambda_1=0.\\
\end{array}\right.
\end{align}
We notice that if $\alpha+2=2\beta,$~~Eq.(3.37) would reduce to~
\begin{align}
\left\{\begin{array}{l}
(\beta-1)^2\lambda_2+(2-\beta)(\beta-1)^2\lambda_3=0,\\
\\
(\beta-1)^3\lambda_2=0,\\
\\
(\beta-1)\lambda_1=0.\\
\end{array}\right.
\end{align}
So when $\beta=1,$~~then $\lambda_1=0,\lambda_2\in R,\lambda_3\in R$.~~When $\beta\neq1,\beta\neq2,$~~then $\lambda_1=\lambda_2=\lambda_3=0$.~~When~~$\beta\neq1,\beta=2,$~~then $\lambda_1=\lambda_2=0,\lambda_3\in R$.~~\\
But if $\alpha+2\neq2\beta,$~~then we have $\lambda_1=0$.\\
We assume that
\begin{align}
\left\{\begin{array}{l}
A=\left[(\frac{\alpha}{2}+1+\beta)(\frac{1}{2}\beta-1)+1\right],\\
\\
B=\{(\alpha-\beta)\left[(2-\beta)(\frac{\alpha}{2}+1)-1\right]+\frac{1}{2}(\frac{\alpha}{2}+1-\beta)\},\\
\\
C=\{\alpha[(\beta-2)(\frac{\alpha}{2}+1)+1]+\frac{1}{2}(\beta-\frac{\alpha}{2}-1)\},\\
\\
D=\frac{1}{2}\beta(\beta-\frac{\alpha}{2}-1).\\
\end{array}\right.
\end{align}
We get
\begin{align}
\left|\begin{array}{cc}
A&B\\
C&D\\
\end{array}\right|
=4\alpha^4\beta^2-16\alpha^4\beta+16\alpha^4-4\alpha^3\beta^3+32\alpha^3\beta^2-68\alpha^3\beta+40\alpha^3-16\alpha^2\beta^3+73\alpha^2\beta^2-\\\notag
98\alpha^2\beta+41\alpha^2-20\alpha\beta^3+60\alpha\beta^2-60\alpha\beta+20\alpha+4\beta^4-16\beta^3+24\beta^2-16\beta+4\notag
\end{align}
Next if \begin{align}
\left|\begin{array}{cc}
A&B\\
C&D\\
\end{array}\right|
=0,
\end{align}
we get $\lambda_3=-\frac{C}{D}\lambda_2$.\\
We continue to assume that $\eta=1,$~~we have
\begin{align}
\left\{\begin{array}{l}
\left[(\frac{\alpha}{2}+1-\beta)(\frac{1}{2}\beta+1)+1\right]\lambda_2+\{(\alpha-\beta)\left[(-2-\beta)(\frac{\alpha}{2}-1)-1\right]+\frac{1}{2}(\frac{\alpha}{2}-1-\beta)\}\lambda_3=0,\\
\\
\{\alpha[(\beta+2)(\frac{\alpha}{2}+1)-1]+\frac{1}{2}(\beta-\frac{\alpha}{2}+1)\}\lambda_2+\frac{1}{2}\beta(\beta-\frac{\alpha}{2}+1)\lambda_3=0,\\
\\
\left[(-2-\beta)(\frac{\alpha}{2}-1+\beta)-2\right]\lambda_1=0,\\
\\
\left[(-2-\beta)(\frac{\alpha}{2}-1)-1\right]\lambda_1=0,\\
\\
(\frac{\alpha}{2}-1-\beta)\lambda_1=0.\\
\end{array}\right.
\end{align}
We notice that if $\alpha=2+2\beta,$~~Eq.(3.42) would reduce to~
\begin{align}
\left\{\begin{array}{l}
(\beta+1)^2\lambda_2-(2+\beta)(\beta-1)^2\lambda_3=0,\\
\\
(\beta+1)^3\lambda_2=0,\\
\\
(\beta+1)\lambda_1=0.\\
\end{array}\right.
\end{align}
So when~~ $\beta=-1,$~~then $\lambda_1=0,\lambda_2\in R,\lambda_3\in R$.~~When~~ $\beta\neq-1,\beta\neq-2,$~~then $\lambda_1=\lambda_2=\lambda_3=0$.~~When$\beta\neq-1,\beta=-2,$~~then $\lambda_1=\lambda_2=0,\lambda_3\in R$.~~\\
But if $\alpha+2\neq2\beta,$~~then we have $\lambda_1=0$.\\
We assume that
\begin{align}
\left\{\begin{array}{l}
A=\left[(\frac{\alpha}{2}-1+\beta)(\frac{1}{2}\beta+1)+1\right],\\
\\
B=\{(\alpha-\beta)\left[(-2-\beta)(\frac{\alpha}{2}-1)-1\right]+\frac{1}{2}(\frac{\alpha}{2}-1-\beta)\},\\
\\
C=\{\alpha[(\beta+2)(\frac{\alpha}{2}-1)+1]+\frac{1}{2}(\beta-\frac{\alpha}{2}+1)\},\\
\\
D=\frac{1}{2}\beta(\beta-\frac{\alpha}{2}+1).\\
\end{array}\right.
\end{align}
We get
\begin{align}
\left|\begin{array}{cc}
A&B\\
C&D\\
\end{array}\right|
=4\alpha^4\beta^2+16\alpha^4\beta+16\alpha^4-4\alpha^3\beta^3-32\alpha^3\beta^2-68\alpha^3\beta-40\alpha^3+16\alpha^2\beta^3+73\alpha^2\beta^2+\\\notag
98\alpha^2\beta+41\alpha^2-20\alpha\beta^3-60\alpha\beta^2-60\alpha\beta-20\alpha+4\beta^4+16\beta^3+24\beta^2+16\beta+4\notag
\end{align}
Next if \begin{align}
\left|\begin{array}{cc}
A&B\\
C&D\\
\end{array}\right|
=0,
\end{align}
we get $\lambda_3=-\frac{C}{D}\lambda_2$.\\
\end{pf}
By (2.89) in \cite{C2}, we have
\vskip 0.5 true cm
\begin{lem}
Ricci symmetric tensors of $(G_4,g,J)$ ~~associated to Kobayashi-Nomizu connection are given by
\begin{align}
\widetilde{{\rm Ric}}^1\left(\begin{array}{c}
e_1\\
e_2\\
e_3
\end{array}\right)=\left(\begin{array}{ccc}
-[1+(\beta-2\eta)(b_3-b_1)]&0&0\\
0&-[1+(\beta-2\eta)(b_2+b_3)]&\frac{b_1+\beta-\alpha-b_3}{2}\\
0&\frac{\alpha+b_3-b_1-\beta}{2}&0
\end{array}\right)\left(\begin{array}{c}
e_1\\
e_2\\
e_3
\end{array}\right).
\end{align}
\end{lem}
By (3.47) and (2.3), we have
\begin{lem}
For $G_4$, the following equalities hold
\begin{align}
&{\rm L}_{e_1} \widetilde{\rm Ric}^1(e_1,e_2)=0,~~{\rm L}_{e_1} \widetilde{\rm Ric}^1(e_1,e_3)=0~~,{\rm L}_{e_1} \widetilde{\rm Ric}^1(e_2,e_2)=-2-\alpha(\beta-2\eta),\\\notag
&{\rm L}_{e_1} \widetilde{\rm Ric}^1(e_2,e_3)=-\beta[1+\alpha(\beta-2\eta)]~~,{\rm L}_{e_1} \widetilde{\rm Ric}^1(e_3,e_3)=\alpha\beta,~~{\rm L}_{e_2} \widetilde{\rm Ric}^1(e_1,e_1)=0~~\\\notag
&{\rm L}_{e_2} \widetilde{\rm Ric}^1(e_1,e_2)=1+\frac{1}{2}\alpha(\beta-2\eta),~~{\rm L}_{e_2} \widetilde{\rm Ric}^1(e_1,e_3)=\alpha[\frac{1}{2}+\beta(\beta-2\eta)]~~,{\rm L}_{e_2} \widetilde{\rm Ric}^1(e_2,e_3)=0,\\\notag
&{\rm L}_{e_2} \widetilde{\rm Ric}^1(e_3,e_3)=0~~,{\rm L}_{e_3} \widetilde{\rm Ric}^1(e_1,e_1)=0,~~{\rm L}_{e_3} \widetilde{\rm Ric}^1(e_1,e_2)=\beta-\frac{\alpha}{2}~~\\\notag
&{\rm L}_{e_3} \widetilde{\rm Ric}^1(e_1,e_3)=-\frac{1}{2}\alpha\beta,~~{\rm L}_{e_3} \widetilde{\rm Ric}^1(e_2,e_2)=0~~,{\rm L}_{e_3} \widetilde{\rm Ric}^1(e_2,e_3)=0.~~\\\notag
\end{align}
\end{lem}
By (3.48) and definition 2.3, we have
\vskip 0.5 true cm

If the unimodular Lorentzian Lie group $(G_4,g,J)$ admits left-invariant Ricci collineations associated to the Kobayashi-Nomizu connection $\nabla^1$, then
${\rm L}_{\xi} \widetilde{\rm Ric}_{ij}^1=0,i=1,2,3~~ and ~~j=1,2,3$, so
\begin{align}
\left\{\begin{array}{l}
\left[1+\frac{1}{2}\alpha(\beta-2\eta)\right]\lambda_2+(\beta-\frac{\alpha}{2})\lambda_3=0,\\
\\
\alpha\left[\frac{1}{2}+\beta(\beta-2\eta)\right]\lambda_2-\frac{1}{2}\alpha\beta\lambda_3=0,\\
\\
\left[-2-\alpha(\beta-2\eta)\right]\lambda_1=0,\\
\\
\beta\left[1+\alpha(\beta-2\eta)\right]\lambda_1=0,\\
\\
\alpha\beta\lambda_1=0.\\
\end{array}\right.
\end{align}
By (3.49), we get
\vskip 0.5 true cm
\begin{thm}
the unimodular Lorentzian Lie group $(G_4,g,J)$ admits left-invariant Ricci collineations associated to Kobayashi-Nomizu connection $\nabla^1$ if and only if\\
$(1)\alpha=\beta=\lambda_1=\lambda_2=0,\mathbb{V}_{\mathbb{R}\mathbb{C}}=\langle e_3\rangle,$\\
\\
$(2)\alpha=\lambda_1=0,\beta\neq0,\mathbb{V}_{\mathbb{R}\mathbb{C}}=\langle e_2-\frac{1}{\beta} e_3\rangle,$\\
\\
$(3)\alpha\neq0,\beta=0,\alpha\eta=1,\mathbb{V}_{\mathbb{R}\mathbb{C}}=\langle e_1\rangle,$\\
\\
$(4)\alpha\neq0,\beta\neq0,\eta=1,\alpha=4\beta,\mathbb{V}_{\mathbb{R}\mathbb{C}}=\langle e_2+\frac{\alpha^2-8\alpha+8}{2\alpha}e_3\rangle,$\\
\\
$(5)\alpha\neq0,\beta\neq0,\eta=1,\alpha=2,\beta=1,\lambda_2=-\lambda_3,\mathbb{V}_{\mathbb{R}\mathbb{C}}=\langle e_2-e_3\rangle,$\\
\\
$(6)\alpha\neq0,\beta\neq0,\eta=-1,\alpha=4\beta,\mathbb{V}_{\mathbb{R}\mathbb{C}}=\langle e_2+\frac{\alpha^2+8\alpha+8}{2\alpha}e_3\rangle,$\\
\\
$(7)\alpha\neq0,\beta\neq0,\eta=-1,\alpha=-2,\beta=-1,\lambda_2=\lambda_3,\mathbb{V}_{\mathbb{R}\mathbb{C}}=\langle e_2+e_3\rangle,$\\
\end{thm}
\begin{pf}
  We analyze each one of these factors by separate and because $b_1=\frac{\alpha}{2}+\eta-\beta,b_3=\frac{\alpha}{2}+\eta,$~~if $\alpha=0$,~~then Eq.(3.49) would reduce to
\begin{align}
\left\{\begin{array}{l}
\lambda_2+\beta\lambda_3=0,\\
\\
\lambda_1=0,\\
\end{array}\right.
\end{align}
 Furthermore if $\beta=0,$~~then we have $\lambda_1=\lambda_2=0,$~~and if $\beta\neq0,$~~then we have $ \lambda_3=-\frac{1}{\beta}\lambda_2.$~~\\
 But if $\alpha\neq0,\beta=0,$~~we get
 \begin{align}
\left\{\begin{array}{l}
(-2+2\alpha\eta)\lambda_1=0,\\
\\
(1-\alpha\eta)\lambda_2-\frac{1}{2}\alpha\lambda_3=0,\\
\\
\lambda_2=0.\\
\end{array}\right.
\end{align}
as long as $\alpha\eta=1,$~~then $\lambda_1\in R,\mathbb{V}_{\mathbb{R}\mathbb{C}}=\langle e_1\rangle$,~~otherwise $\mathbb{V}_{\mathbb{R}\mathbb{C}}$~~is empty set.~~\\
Considering $\alpha\neq0,\beta\neq0,$~~then naturely $\lambda_1=0,$~~and when $\eta=1,$~~we have
\begin{align}
\left\{\begin{array}{l}
\left[1+\frac{1}{2}\alpha(\beta-2)\right]\lambda_2+(\beta-\frac{\alpha}{2})\lambda_3=0,\\
\\
\left[\frac{1}{2}+\beta(\beta-2)\right]\lambda_2-\frac{1}{2}\beta\lambda_3=0,\\
\end{array}\right.
\end{align}
We assume that
\begin{align}
\left\{\begin{array}{l}
A=\left[1+\frac{1}{2}\alpha(\beta-2)\right],\\
\\
B=\beta-\frac{\alpha}{2},\\
\\
C=\left[\frac{1}{2}+\beta(\beta-2)\right],\\
\\
D=-\frac{1}{2}\beta.\\
\end{array}\right.
\end{align}
We get
\begin{align}
\left|\begin{array}{cc}
A&B\\
C&D\\
\end{array}\right|
=(\frac{1}{4}\alpha-\beta)(\beta-1)^2
\end{align}
If $\beta=1,$~~we have~~$(1-\frac{1}{2}\alpha)(\lambda_2+\lambda_3)=0,$~~so $\lambda_2+\lambda_3 ~~equals~~ to ~~0 $~~with $\alpha=2.$~~\\
If $\alpha=4\beta,$~~we have~~$(1+2\beta^2-4\beta)-\beta\lambda_3=0,$~~so $\lambda_3=\frac{1+2\beta^2-4\beta}{\beta}\lambda_2=\frac{\alpha^2-8\alpha+8}{2\alpha}\lambda_2.$\\
Now if $\eta=-1,$~~we have
\begin{align}
\left\{\begin{array}{l}
\left[1+\frac{1}{2}\alpha(\beta+2)\right]\lambda_2+(\beta-\frac{\alpha}{2})\lambda_3=0,\\
\\
\left[\frac{1}{2}+\beta(\beta+2)\right]\lambda_2-\frac{1}{2}\beta\lambda_3=0,\\
\end{array}\right.
\end{align}
We assume that
\begin{align}
\left\{\begin{array}{l}
A=\left[1+\frac{1}{2}\alpha(\beta+2)\right],\\
\\
B=\beta-\frac{\alpha}{2},\\
\\
C=\left[\frac{1}{2}+\beta(\beta+2)\right],\\
\\
D=-\frac{1}{2}\beta.\\
\end{array}\right.
\end{align}
We get
\begin{align}
\left|\begin{array}{cc}
A&B\\
C&D\\
\end{array}\right|
=(\frac{1}{4}\alpha-\beta)(\beta+1)^2
\end{align}
If $\beta=-1,$~~we have$(1+\frac{1}{2}\alpha)(\lambda_2-\lambda_3)=0,$~~so $\lambda_2+\lambda_3~~ equals~~ to~~ 0 $~~with ~$\alpha=-2.$~~\\
If $\alpha=4\beta,$~~we have$(1+2\beta^2+4\beta)-\beta\lambda_3=0,$~~so $\lambda_3=\frac{1+2\beta^2+4\beta}{\beta}\lambda_2=\frac{\alpha^2+8\alpha+8}{2\alpha}\lambda_3.$\\

In summary,the unimodular Lorentzian Lie group $(G_4,g,J)$ admits left-invariant Ricci collineations associated to Kobayashi-Nomizu connection $\nabla^1$ if and only if\\
$(1)\alpha=\beta=\lambda_1=\lambda_2=0,\mathbb{V}_{\mathbb{R}\mathbb{C}}=\langle e_3\rangle,$\\
\\
$(2)\alpha=\lambda_1=0,\beta\neq0,\mathbb{V}_{\mathbb{R}\mathbb{C}}=\langle e_2-\frac{1}{\beta} e_3\rangle,$\\
\\
$(3)\alpha\neq0,\beta=0,\alpha\eta=1,\mathbb{V}_{\mathbb{R}\mathbb{C}}=\langle e_1\rangle,$\\
\\
$(4)\alpha\neq0,\beta\neq0,\eta=1,\alpha=4\beta,\mathbb{V}_{\mathbb{R}\mathbb{C}}=\langle e_2+\frac{\alpha^2-8\alpha+8}{2\alpha}e_3\rangle,$\\
\\
$(5)\alpha\neq0,\beta\neq0,\eta=1,\alpha=2,\beta=1,\lambda_2=-\lambda_3,\mathbb{V}_{\mathbb{R}\mathbb{C}}=\langle e_2-e_3\rangle,$\\
\\
$(6)\alpha\neq0,\beta\neq0,\eta=-1,\alpha=4\beta,\mathbb{V}_{\mathbb{R}\mathbb{C}}=\langle e_2+\frac{\alpha^2+8\alpha+8}{2\alpha}e_3\rangle,$\\
\\
$(7)\alpha\neq0,\beta\neq0,\eta=-1,\alpha=-2,\beta=-1,\lambda_2=\lambda_3,\mathbb{V}_{\mathbb{R}\mathbb{C}}=\langle e_2+e_3\rangle,$\\
\end{pf}

\vskip 1 true cm

\section{ Left-invariant Ricci collineations associated to canonical connections and Kobayashi-Nomizu connections on three-dimensional non-unimodular Lorentzian Lie groups}

\vskip 0.5 true cm
\noindent{\bf 4.1 Left-invariant Ricci collineation of $G_5$}\\
\vskip 0.5 true cm
By \cite{A1}, we have for $G_5$, there exists a pseudo-orthonormal basis $\{e_1,e_2,e_3\}$ with $e_3$ timelike such that the Lie
algebra of $G_5$ satisfies
\begin{equation}
[e_1,e_2]=0,~~[e_1,e_3]=\alpha e_1+\beta e_2,~~[e_2,e_3]=\gamma e_1+\delta e_2,\alpha+\delta\neq0,\alpha\gamma+\beta\delta=0,~~.\notag
\end{equation}
\vskip 0.5 true cm

By (3.5) in \cite{C2}, we have
\vskip 0.5 true cm
\begin{lem}
Ricci symmetric tensors of $(G_5,g,J)$  associated to canonical connection are given by
\begin{align}
\widetilde{{\rm Ric}}^0\left(\begin{array}{c}
e_1\\
e_2\\
e_3
\end{array}\right)=\left(\begin{array}{ccc}
0&0&0\\
0&0&0\\
0&0&0
\end{array}\right)\left(\begin{array}{c}
e_1\\
e_2\\
e_3
\end{array}\right).
\end{align}
\end{lem}

\begin{lem}
For $G_5$, the following equalities hold~~${\rm L}_{\xi} \widetilde{\rm Ric}^0(i,j)=0,for ~~i=1,2,3~~~and~~~j=1,2,3.$
\end{lem}
\vskip 0.1 true cm
If the unimodular Lorentzian Lie group $(G_5,g,J)$ admits left-invariant Ricci collineations associated to the canonical connection $\nabla^0$, then
${\rm L}_{\xi} \widetilde{\rm Ric}_{ij}^0=0,i=1,2,3~~~and~~~j=1,2,3$, so there are no constraint conditions for $G_5$.\\
\vskip 0.1 true cm
Naturelly, we get
\vskip 0.5 true cm
\begin{thm}
the unimodular Lorentzian Lie group $(G_5,g,J)$ admits definitely left-invariant Ricci collineations associated to the canonical connection $\nabla^0$ .\\
\end{thm}
\indent By \cite{C2}, we have
\vskip 0.5 true cm
\begin{lem}
Ricci symmetric tensors of $(G_5,g,J)$ ~~associated to Kobayashi-Nomizu connection are given by
\begin{align}
\widetilde{{\rm Ric}}^1\left(\begin{array}{c}
e_1\\
e_2\\
e_3
\end{array}\right)=\left(\begin{array}{ccc}
0&0&0\\
0&0&0\\
0&0&0
\end{array}\right)\left(\begin{array}{c}
e_1\\
e_2\\
e_3
\end{array}\right).
\end{align}
\end{lem}
\begin{lem}
For $G_5$, the following equalities hold~~~${\rm L}_{\xi} \widetilde{\rm Ric}^1(i,j)=0,i=1,2,3~~~and~~~j=1,2,3.$
\end{lem}
\vskip 0.1 true cm
If the unimodular Lorentzian Lie group $(G_5,g,J)$ admits left-invariant Ricci collineations associated to Kobayashi-Nomizu connection $\nabla^1$, then
${\rm L}_{\xi} \widetilde{\rm Ric}_{ij}^1=0,i=1,2,3~~ and ~~j=1,2,3$, so there are no constraint conditions for $G_5$.\\
Naturelly, we get
\vskip 0.5 true cm
\begin{thm}
the unimodular Lorentzian Lie group $(G_5,g,J)$ admits definitely left-invariant Ricci collineations associated to Kobayashi-Nomizu connection $\nabla^1$ .\\
\end{thm}
\vskip 0.3 true cm
\noindent{\bf 4.2 Left-invariant Ricci collineation of $G_6$}\\
\vskip 0.5 true cm
By \cite{A1}, we have for $G_6$, there exists a pseudo-orthonormal basis $\{e_1,e_2,e_3\}$ with $e_3$ timelike such that the Lie
algebra of $G_6$ satisfies
\begin{equation}
[e_1,e_2]=\alpha e_2+\beta e_3,~~[e_1,e_3]=\gamma e_2+\delta e_3,~~[e_2,e_3]=0,\alpha+\delta\neq0,\alpha\gamma-\beta\delta=0~~.\notag
\end{equation}
\vskip 0.5 true cm

By (3.18) in \cite{C2}, we have
\vskip 0.5 true cm
\begin{lem}
Ricci symmetric tensors of $(G_6,g,J)$  associated to canonical connection are given by
\begin{align}
\widetilde{{\rm Ric}}^0\left(\begin{array}{c}
e_1\\
e_2\\
e_3
\end{array}\right)=\left(\begin{array}{ccc}
\frac{1}{2}\beta(\beta-\gamma)-\alpha^2&0&0\\
0&\frac{1}{2}\beta(\beta-\gamma)-\alpha^2&\frac{1}{2}[\gamma\alpha-\frac{1}{2}\delta(\beta-\gamma)]\\
0&\frac{1}{2}[-\gamma\alpha+\frac{1}{2}\delta(\beta-\gamma)]&0
\end{array}\right)\left(\begin{array}{c}
e_1\\
e_2\\
e_3
\end{array}\right).
\end{align}
\end{lem}
By (4.3) and (2.2)
\begin{lem}
For $G_6$, the following equalities hold
\vskip 0.1 true cm
\begin{align}
&{\rm L}_{e_1} \widetilde{\rm Ric}^0(e_1,e_2)=0,~~{\rm L}_{e_1} \widetilde{\rm Ric}^0(e_1,e_3)=0~~,{\rm L}_{e_1} \widetilde{\rm Ric}^0(e_2,e_2)=2\alpha^3+\alpha\beta\gamma-\beta(\beta-\gamma)(\alpha+\frac{1}{2}\delta),\\\notag
&{\rm L}_{e_1} \widetilde{\rm Ric}^0(e_2,e_3)=\frac{1}{2}\alpha\gamma(3\alpha+\delta)+\frac{1}{2}(\gamma-\beta)(\alpha\delta+\delta^2+\beta\gamma),~~{\rm L}_{e_1} \widetilde{\rm Ric}^0(e_3,e_3)=\gamma\left[\alpha\gamma+\frac{1}{2}\delta(\gamma-\beta)\right],~~\\\notag
&{\rm L}_{e_2} \widetilde{\rm Ric}^0(e_1,e_1)=0~~,{\rm L}_{e_2} \widetilde{\rm Ric}^0(e_1,e_2)=-\alpha^3-\frac{1}{2}\alpha\beta\gamma+\frac{1}{2}\beta(\beta-\gamma)(\alpha+\frac{1}{2}\delta),~~\\\notag
&{\rm L}_{e_2} \widetilde{\rm Ric}^0(e_1,e_3)=\frac{1}{2}\alpha\left[-\alpha\gamma+\frac{1}{2}\delta(\beta-\gamma)\right]~~,{\rm L}_{e_2} \widetilde{\rm Ric}^0(e_2,e_3)=0,\\\notag
&{\rm L}_{e_2} \widetilde{\rm Ric}^0(e_3,e_3)=0~~,{\rm L}_{e_3} \widetilde{\rm Ric}^0(e_1,e_1)=0,~~{\rm L}_{e_3} \widetilde{\rm Ric}^0(e_1,e_2)=-\alpha^2\gamma-\frac{1}{2}\alpha\gamma\delta+\frac{1}{2}(\beta-\gamma)(\beta\gamma+\frac{1}{2}\delta^2)~~\\\notag
&{\rm L}_{e_3} \widetilde{\rm Ric}^0(e_1,e_3)=\frac{1}{2}\gamma\left[-\alpha\gamma+\frac{1}{2}\delta(\beta-\gamma)\right],~~{\rm L}_{e_3} \widetilde{\rm Ric}^0(e_2,e_2)=0~~,{\rm L}_{e_3} \widetilde{\rm Ric}^0(e_2,e_3)=0.~~\\\notag
\end{align}
\end{lem}
By (4.4) and definition 2.3, we have
\vskip 0.5 true cm

If the unimodular Lorentzian Lie group $(G_6,g,J)$ is left-invariant Ricci collineations associated to the canonical connection $\nabla^0$, then
${\rm L}_{\xi} \widetilde{\rm Ric}_{ij}^0=0,i=1,2,3~~ and~~ j=1,2,3$, so
\\
\begin{align}
\left\{\begin{array}{l}
\left[-\alpha^3-\frac{1}{2}\alpha\beta\gamma+\frac{1}{2}\beta(\beta-\gamma)(\alpha+\frac{1}{2}\delta)\right]\lambda_2+\left[-\alpha^2\gamma-\frac{1}{2}\alpha\gamma\delta+\frac{1}{2}(\beta-\gamma)(\beta\gamma+\frac{1}{2}\delta^2)\right]\lambda_3=0,\\
\\
\frac{1}{2}\alpha\left[-\alpha\gamma+\frac{1}{2}\delta(\beta-\gamma)\right]\lambda_2+\frac{1}{2}\gamma\left[-\alpha\gamma+\frac{1}{2}\delta(\beta-\gamma)\right]\lambda_3=0,\\
\\
\left[2\alpha^3+\alpha\beta\gamma-\beta(\beta-\gamma)(\alpha+\frac{1}{2}\delta)\right]\lambda_1=0,\\
\\
\left[\frac{1}{2}\alpha\gamma(3\alpha+\delta)+\frac{1}{2}(\gamma-\beta)(\alpha\delta+\delta^2+\beta\gamma)\right]\lambda_1=0,\\
\\
\gamma\left[\alpha\gamma+\frac{1}{2}\delta(\gamma-\beta)\right]\lambda_1=0.\\
\end{array}\right.
\end{align}
By solving (4.5), we get
\vskip 0.5 true cm
\begin{thm}
the unimodular Lorentzian Lie group $(G_6,g,J)$ is left-invariant Ricci collineations associated to the canonical connection $\nabla^0$ if and only if\\
\\
$(1)\gamma=0,\alpha(2\alpha^2-\beta^2)=0,\mathbb{V}_{\mathbb{R}\mathbb{C}}=\langle e_1,e_2,e_3\rangle,$\\
\\
$(2)\gamma=0,\alpha(2\alpha^2-\beta^2)\neq0,,\mathbb{V}_{\mathbb{R}\mathbb{C}}=\langle e_3\rangle,$\\
\\
$(3)\gamma\neq0,\alpha=\beta=0,\delta\neq0,\lambda_1=\lambda_3=0,\mathbb{V}_{\mathbb{R}\mathbb{C}}=\langle e_2\rangle,$\\
\\
$(4)\gamma\neq0,\alpha\neq0,\alpha+\beta=0,\gamma+\delta=0,\delta\neq0,\lambda_1=0,\lambda_3=\frac{\alpha}{\delta}\lambda_2,\mathbb{V}_{\mathbb{R}\mathbb{C}}=\langle e_2+\frac{\alpha}{\delta} e_3\rangle,$\\
\\
$(5)\gamma\neq0,\alpha\neq0,\alpha=\beta,\gamma=\delta,\delta\neq0,\lambda_1=0,\lambda_3=-\frac{\alpha}{\delta}\lambda_2,\mathbb{V}_{\mathbb{R}\mathbb{C}}=\langle e_2-\frac{\alpha}{\delta} e_3\rangle.$\\
\end{thm}
\begin{pf}
  We analyze each one of these factors by separate ,if $\gamma\neq0,\alpha=0$,~~we have
\begin{align}
\left\{\begin{array}{l}
\frac{1}{2}\beta\delta(\beta-\gamma)\lambda_2+(\beta-\gamma)(\beta\gamma+\frac{1}{2}\gamma^2)\lambda_3=0,\delta(\beta-\gamma)\lambda_3=0,\\
\\
\beta\delta(\beta-\gamma)\lambda_1=0,(\beta-\gamma)(\beta\gamma+\gamma^2)\lambda_1=0,\\
\\
\delta(\beta-\gamma)\lambda_1=0,\alpha\gamma-\beta\delta=0,\alpha+\delta\neq0.\\
\end{array}\right.
\end{align}
Notice that $\delta\neq0$~~for $\alpha+\delta\neq0$~~and furthermore $\beta=0$~~for $\alpha\gamma-\beta\delta=0$~~,so Eq.(4.6) would reduce to
\begin{align}
\left\{\begin{array}{l}
(\beta\gamma+\frac{1}{2}\gamma^2)\lambda_3=0,\\
\\
\delta\gamma\lambda_3=0,\\
\\
\beta\gamma\lambda_1=0.\\
\end{array}\right.
\end{align}
Naturelly we get~~$\lambda_1=\lambda_3=0,\lambda_2 \in R$.\\
If $\gamma\neq0,\alpha\neq0,$~~we have $\beta\neq0,\delta\neq0,\left[\alpha\gamma+\frac{1}{2}\delta(\gamma-\beta)\right]\neq0,$~~otherwise $\alpha\gamma-\beta\delta\neq0.$~~\\
Naturelly we have $\lambda_1=0$~~and we assume that
\begin{align}
\left\{\begin{array}{l}
A=-\alpha^3-\frac{1}{2}\alpha\beta\gamma+\frac{1}{2}\beta(\beta-\gamma)(\alpha+\frac{1}{2}\delta),\\
\\
B=-\alpha^2\gamma-\frac{1}{2}\alpha\gamma\delta+\frac{1}{2}(\beta-\gamma)(\beta\gamma+\frac{1}{2}\delta^2),\\
\\
C=\frac{1}{2}\alpha\left[-\alpha\gamma+\frac{1}{2}\delta(\beta-\gamma)\right],\\
\\
D=\frac{1}{2}\gamma\left[-\alpha\gamma+\frac{1}{2}\delta(\beta-\gamma)\right].\\
\end{array}\right.
\end{align}
We get
\begin{align}
\left|\begin{array}{cc}
A&B\\
C&D\\
\end{array}\right|
=\frac{1}{4}\gamma(\alpha+\delta)(\alpha\delta-\beta\gamma)=0,
\end{align}
solving
\begin{align}
\left\{\begin{array}{l}
\alpha\gamma-\beta\delta=0,\\
\\
\alpha\delta-\beta\gamma=0.\\
\\
\end{array}\right.
\end{align}
we get~~$\alpha+\beta=0,\gamma+\delta=0,\lambda_3=\frac{\alpha}{\delta}\lambda_2 ~~or~~  \alpha=\beta,\gamma=\delta, \lambda_3=-\frac{\alpha}{\delta}\lambda_2.$~~\\
If $\gamma=0,$~~we naturely have $\beta\delta=0,$~~ then we get
\begin{align}
\left\{\begin{array}{l}
(2\alpha^3-\alpha\beta^2)\lambda_1=0,\\
\\
(2\alpha^3-\alpha\beta^2)\lambda_2=0.\\
\\
\end{array}\right.
\end{align}
When $\alpha(2\alpha^2-\beta^2)=0,$~~we get $\lambda_1\in R,\lambda_2\in R,\lambda_3\in R$~~and when $\alpha(2\alpha^2-\beta^2)\neq0,$~~we get $\lambda_1=\lambda_2=0,\lambda_3\in R.$~~\\
 In summary,the unimodular Lorentzian Lie group $(G_6,g,J)$ is left-invariant Ricci collineations associated to the canonical connection $\nabla^0$ if and only if\\
\\
$(1)\gamma=0,\alpha(2\alpha^2-\beta^2)=0,\mathbb{V}_{\mathbb{R}\mathbb{C}}=\langle e_1,e_2,e_3\rangle,$\\
\\
$(2)\gamma=0,\alpha(2\alpha^2-\beta^2)\neq0,,\mathbb{V}_{\mathbb{R}\mathbb{C}}=\langle e_3\rangle,$\\
\\
$(3)\gamma\neq0,\alpha=\beta=0,\delta\neq0,\lambda_1=\lambda_3=0,\mathbb{V}_{\mathbb{R}\mathbb{C}}=\langle e_2\rangle,$\\
\\
$(4)\gamma\neq0,\alpha\neq0,\alpha+\beta=0,\gamma+\delta=0,\delta\neq0,\lambda_1=0,\lambda_3=\frac{\alpha}{\delta}\lambda_2,\mathbb{V}_{\mathbb{R}\mathbb{C}}=\langle e_2+\frac{\alpha}{\delta} e_3\rangle,$\\
\\
$(5)\gamma\neq0,\alpha\neq0,\alpha=\beta,\gamma=\delta,\delta\neq0,\lambda_1=0,\lambda_3=-\frac{\alpha}{\delta}\lambda_2,\mathbb{V}_{\mathbb{R}\mathbb{C}}=\langle e_2-\frac{\alpha}{\delta} e_3\rangle.$\\
\end{pf}
\indent By (3.23) in \cite{C2}, we have
\vskip 0.5 true cm
\begin{lem}
Ricci symmetric tensors of $(G_6,g,J)$ ~~associated to Kobayashi-Nomizu connection are given by
\begin{align}
\widetilde{{\rm Ric}}^1\left(\begin{array}{c}
e_1\\
e_2\\
e_3
\end{array}\right)=\left(\begin{array}{ccc}
-(\alpha^2+\beta\gamma)&0&0\\
0&-\alpha^2&0\\
0&0&0
\end{array}\right)\left(\begin{array}{c}
e_1\\
e_2\\
e_3
\end{array}\right).
\end{align}
\end{lem}
By(4.12) and (2.3), we have
\begin{lem}
For $G_6$, the following equalities hold
\begin{align}
&{\rm L}_{e_1} \widetilde{\rm Ric}^1(e_1,e_2)=0,~~{\rm L}_{e_1} \widetilde{\rm Ric}^1(e_1,e_3)=0~~,{\rm L}_{e_1} \widetilde{\rm Ric}^1(e_2,e_2)=2\alpha^3,\\\notag
&{\rm L}_{e_1} \widetilde{\rm Ric}^1(e_2,e_3)=\alpha^2\gamma,~~{\rm L}_{e_1} \widetilde{\rm Ric}^1(e_3,e_3)=0,~~{\rm L}_{e_2} \widetilde{\rm Ric}^1(e_1,e_1)=0,~~\\\notag
&{\rm L}_{e_2} \widetilde{\rm Ric}^1(e_1,e_2)=-\alpha^3,~~{\rm L}_{e_2} \widetilde{\rm Ric}^1(e_1,e_3)=0~~,{\rm L}_{e_2} \widetilde{\rm Ric}^1(e_2,e_3)=0,\\\notag
&{\rm L}_{e_2} \widetilde{\rm Ric}^1(e_3,e_3)=0~~,{\rm L}_{e_3} \widetilde{\rm Ric}^1(e_1,e_1)=0,~~{\rm L}_{e_3} \widetilde{\rm Ric}^1(e_1,e_2)=-\alpha^2\gamma,~~\\\notag
&{\rm L}_{e_3} \widetilde{\rm Ric}^1(e_1,e_3)=0,~~{\rm L}_{e_3} \widetilde{\rm Ric}^1(e_2,e_2)=0~~,{\rm L}_{e_3} \widetilde{\rm Ric}^1(e_2,e_3)=0.~~\\\notag
\end{align}
\end{lem}
By (4.13) and definition 2.3, we have
\vskip 0.5 true cm
If the unimodular Lorentzian Lie group $(G_6,g,J)$ admits left-invariant Ricci collineations associated to Kobayashi-Nomizu connection $\nabla^1$, then
${\rm L}_{\xi} \widetilde{\rm Ric}_{ij}^1=0,for~~i=1,2,3~~ and ~~j=1,2,3$,~~ so
\begin{align}
\left\{\begin{array}{l}
\alpha^3\lambda_2+\alpha^2\gamma\lambda_3=0,\\
\\
2\alpha^3\lambda_1=0,\\
\\
\alpha^2\gamma\lambda_1=0.\\
\end{array}\right.
\end{align}
By solving (4.14), we get
\vskip 0.5 true cm
\begin{thm}
the unimodular Lorentzian Lie group $(G_6,g,J)$ admits left-invariant Ricci collineations associated to Kobayshi-Nomizu connection $\nabla^1$ if and only if\\
$(1)\alpha=\beta=0,\delta\neq0,\mathbb{V}_{\mathbb{R}\mathbb{C}}=\langle e_1,e_2, e_3\rangle,$\\
\\
$(2)\alpha\neq0,\lambda_2=-\frac{\gamma}{\alpha}\lambda_3,\mathbb{V}_{\mathbb{R}\mathbb{C}}=\langle -\frac{\gamma}{\alpha}e_2+e_3\rangle.$\\
\end{thm}
\begin{pf}
  We analyze each one of these factors by separate. If $\alpha=0$,~~we have~~$\delta\neq0,\beta=0$~~for $\alpha+\delta\neq0,\alpha\gamma-\beta\delta=0,$~~and naturely $\lambda_1\in R,\lambda_2\in R,\lambda_3\in R.$~~\\
And if $\alpha\neq0,$~~we get
\begin{align}
\left\{\begin{array}{l}
\alpha\lambda_2+\gamma\lambda_3=0,\\
\\
\lambda_1=0.\\
\end{array}\right.
\end{align}
So $\lambda_2=-\frac{\gamma}{\alpha}\lambda_3, \mathbb{V}_{\mathbb{R}\mathbb{C}}=\langle -\frac{\gamma}{\alpha}e_2+e_3\rangle.$\\
\end{pf}
\vskip 0.5 true cm

\vskip 0.5 true cm
\noindent{\bf 4.3 Left-invariant Ricci collineation of $G_7$}\\
\vskip 0.5 true cm
By \cite{A1}, we have for $G_7$, there exists a pseudo-orthonormal basis $\{e_1,e_2,e_3\}$ with $e_3$ timelike such that the Lie
algebra of $G_7$ satisfies
\begin{equation}
[e_1,e_2]=-\alpha e_1-\beta e_2-\beta e_3,~~[e_1,e_3]=\alpha e_1+\beta e_2+\beta e_3,~~[e_2,e_3]=\gamma e_1+\delta e_2+\delta e_3,~~\alpha+\delta\neq0,~~\alpha\gamma=0.\notag
\end{equation}
\vskip 0.5 true cm

By (3.34) in \cite{C2}, we have
\vskip 0.5 true cm
\begin{lem}
Ricci symmetric tensors of $(G_7,g,J)$  associated to canonical connection are given by
\begin{align}
\widetilde{{\rm Ric}}^0\left(\begin{array}{c}
e_1\\
e_2\\
e_3
\end{array}\right)=\left(\begin{array}{ccc}
-(\alpha^2+\frac{\beta\gamma}{2})&0&\frac{1}{2}(\gamma\alpha+\frac{\delta\gamma}{2})\\
0&-(\alpha^2+\frac{\beta\gamma}{2})&-\frac{1}{2}(\alpha^2+\frac{\beta\gamma}{2})\\
-\frac{1}{2}(\gamma\alpha+\frac{\delta\gamma}{2})&\frac{1}{2}(\alpha^2+\frac{\beta\gamma}{2})&0
\end{array}\right)\left(\begin{array}{c}
e_1\\
e_2\\
e_3
\end{array}\right).
\end{align}
\end{lem}
By (4.16) and (2.2),we have
\begin{lem}
For $G_7$, the following equalities hold
\vskip 0.1 true cm
\begin{align}
&{\rm L}_{e_1} \widetilde{\rm Ric}^0(e_1,e_2)=-\alpha^3-\frac{1}{4}\beta\delta\gamma,{\rm L}_{e_1} \widetilde{\rm Ric}^0(e_1,e_3)=\alpha^3+\frac{1}{4}\beta\delta\gamma~~\\\notag
&{\rm L}_{e_1} \widetilde{\rm Ric}^0(e_2,e_2)=-\beta(\alpha^2+\frac{1}{2}\beta\gamma),{\rm L}_{e_1} \widetilde{\rm Ric}^0(e_2,e_3)=\beta(\alpha^2+\frac{1}{2}\beta\gamma),\\\notag
&{\rm L}_{e_1} \widetilde{\rm Ric}^0(e_3,e_3)=-\beta(\alpha^2+\frac{1}{2}\beta\gamma),{\rm L}_{e_2} \widetilde{\rm Ric}^0(e_1,e_1)=2\alpha^3+\frac{1}{2}\beta\delta\gamma~~\\\notag
&{\rm L}_{e_2} \widetilde{\rm Ric}^0(e_1,e_2)=\frac{1}{2}\beta(\alpha^2+\frac{1}{2}\beta\gamma),{\rm L}_{e_2} \widetilde{\rm Ric}^0(e_1,e_3)=\frac{1}{4}\gamma(\delta^2-\beta^2)+\frac{1}{2}\beta(\gamma^2-\alpha^2)~~\\\notag
&{\rm L}_{e_2} \widetilde{\rm Ric}^0(e_2,e_3)=\frac{1}{2}\delta(\alpha^2+\frac{1}{2}\beta\gamma),{\rm L}_{e_2} \widetilde{\rm Ric}^0(e_3,e_3)=\delta(\frac{1}{2}\gamma^2-\alpha^2-\frac{1}{2}\beta\gamma)~~\\\notag
&{\rm L}_{e_3} \widetilde{\rm Ric}^0(e_1,e_1)=-2\alpha^3-\frac{1}{2}\beta\delta\gamma,{\rm L}_{e_3} \widetilde{\rm Ric}^0(e_1,e_2)=-(\gamma+\frac{1}{2}\beta)(\alpha^2+\frac{1}{2}\beta\gamma)-\frac{1}{4}\delta^2\gamma~~\\\notag
&{\rm L}_{e_3} \widetilde{\rm Ric}^0(e_1,e_3)=\frac{1}{2}\beta(\alpha^2+\frac{1}{2}\beta\gamma),{\rm L}_{e_3} \widetilde{\rm Ric}^0(e_2,e_2)=-\delta(\alpha^2+\frac{1}{2}\beta\gamma)~~\\\notag
&{\rm L}_{e_3} \widetilde{\rm Ric}^0(e_2,e_3)=-\frac{1}{4}\delta\gamma^2+\frac{1}{2}\delta(\alpha^2+\frac{1}{2}\beta\gamma).~~\\\notag
\end{align}
\end{lem}
By (4.17) and definition 2.3, we have
\vskip 0.5 true cm

If the unimodular Lorentzian Lie group $(G_7,g,J)$ admits left-invariant Ricci collineations associated to the canonical connection $\nabla^0$, then
${\rm L}_{\xi} \widetilde{\rm Ric}_{ij}^0=0,i=1,2,3~~ and~~ j=1,2,3$, so
\begin{align}
\left\{\begin{array}{l}
(2\alpha^3+\frac{1}{2}\beta\delta\gamma)(\lambda_2-\lambda_3)=0,\\
\\
(-\alpha^3-\frac{1}{4}\beta\delta\gamma)\lambda_1+\frac{1}{2}\beta(\alpha^2+\frac{1}{2}\beta\gamma)\lambda_2+\left[-(\gamma+\frac{1}{2}\beta)(\alpha^2+\frac{1}{2}\beta\gamma)-\frac{1}{4}\delta^2\gamma\right]\lambda_3=0,\\
\\
(\alpha^3+\frac{1}{4}\beta\delta\gamma)\lambda_1+\left[\frac{1}{4}\gamma(\delta^2-\beta^2)+\frac{1}{2}\beta(\gamma^2-\alpha^2)\right]\lambda_2+\frac{1}{2}\beta(\alpha^2+\frac{1}{2}\beta\gamma)\lambda_3=0,\\
\\
\beta(\alpha^2+\frac{1}{2}\beta\gamma)\lambda_1+\delta(\alpha^2+\frac{1}{2}\beta\gamma)\lambda_3=0,\\
\\
\beta(\alpha^2+\frac{1}{2}\beta\gamma)\lambda_1+\frac{1}{2}\delta(\alpha^2+\frac{1}{2}\beta\gamma)\lambda_2+\left[-\frac{1}{4}\delta\gamma^2+\frac{1}{2}\delta(\alpha^2+\frac{1}{2}\beta\gamma)\right]\lambda_3=0,\\
\\
-\beta(\alpha^2+\frac{1}{2}\beta\gamma)\lambda_1+\delta(\frac{1}{2}\gamma^2-\alpha^2-\frac{1}{2}\beta\gamma)\lambda_2=0.\\
\end{array}\right.
\end{align}

By solving (4.18), we get
\vskip 0.5 true cm
\begin{thm}
the unimodular Lorentzian Lie group $(G_7,g,J)$ admits left-invariant Ricci collineations associated to the canonical connection $\nabla^0$ if and only if\\
$(1)\alpha=0,\delta\neq0,\gamma\neq0,\beta=0,\lambda_2=\lambda_3=0,\mathbb{V}_{\mathbb{R}\mathbb{C}}=\langle e_1\rangle,$\\
\\
$(2)\alpha=0,\delta\neq0,\gamma=0,\mathbb{V}_{\mathbb{R}\mathbb{C}}=\langle e_1,e_2,e_3\rangle,$\\
\\
$(3)\alpha\neq0,\gamma=\delta=0,\lambda_1=0,\lambda_2=\lambda_3,\mathbb{V}_{\mathbb{R}\mathbb{C}}=\langle e_2+e_3\rangle.$\\
\end{thm}
\begin{pf}
  We analyze each one of these factors by separate. If $\alpha=0$,~~we have~~$\delta\neq0,$~~for $\alpha+\delta\neq0,\alpha\gamma=0,$~~and
\begin{align}
\left\{\begin{array}{l}
\beta\gamma(\lambda_2-\lambda_3)=0,\\
\\
\beta\delta\gamma\lambda_1-\beta^2\gamma\lambda_2+\left[\beta\gamma(2\gamma+\beta)+\delta^2\gamma\right]\lambda_3=0,\\
\\
\beta\delta\gamma\lambda_1+\left[\gamma(\delta^2-\beta^2)+2\beta\gamma^2\right]\lambda_2+\beta^2\gamma\lambda_3=0,\\
\\
\beta^2\gamma\lambda_1+\beta\gamma\delta\lambda_3=0,\\
\\
2\beta^2\gamma\lambda_1+\beta\gamma\delta\lambda_2+(\beta\gamma\delta-\delta\gamma^2)\lambda_3=0,\\
\\
\beta^2\gamma\lambda_1+(\beta\gamma\delta-\delta\gamma^2)\lambda_2=0.\\
\end{array}\right.
\end{align}
Obviously $\lambda_1\in R,\lambda_2\in R,\lambda_3\in R $~~with $\gamma=0,$~~and when $\beta\neq0$,~~we have $\lambda_2=\lambda_3.$~~Then we have\\
\begin{align}
\left\{\begin{array}{l}
\beta\delta\lambda_1+(\delta^2+2\beta\gamma)\lambda_2=0,\\
\\
\beta^2\lambda_1+\beta\delta\lambda_2=0,\\
\\
\beta^2\lambda_1+(\beta\delta-\delta\gamma)\lambda_2=0.\\
\end{array}\right.
\end{align}
So we get~~$\lambda_1=\lambda_2=\lambda_3=0.$\\
when $\gamma\neq0,$~~Eq.(4.19)would easily reduce  to
\begin{align}
\left\{\begin{array}{l}
\beta(\lambda_2-\lambda_3)=0,\\
\\
\beta\delta\lambda_1-\beta^2\lambda_2+\left[\beta(2\gamma+\beta)+\delta^2\right]\lambda_3=0,\\
\\
\beta\delta\lambda_1+(\delta^2-\beta^2+2\beta\gamma)\lambda_2+\beta^2\lambda_3=0,\\
\\
\beta^2\lambda_1+\beta\delta\lambda_3=0,\\
\\
2\beta^2\lambda_1+\beta\delta\lambda_2+(\beta\delta-\delta\gamma)\lambda_3=0,\\
\\
\beta^2\lambda_1+(\beta\delta-\delta\gamma)\lambda_2=0.\\
\end{array}\right.
\end{align}
simplify this equation, we still get $\lambda_2=\lambda_3=0,\lambda_1\in R$~~with $\beta=0.$~~\\
If $\alpha\neq0,$~~then we have $\gamma=0.$~~And we get
\begin{align}
\left\{\begin{array}{l}
\alpha^3(\lambda_2-\lambda_3)=0,\\
\\
\alpha^3\lambda_1+\frac{1}{2}\alpha^2\beta(\lambda_3-\lambda_2)=0,\\
\\
\beta\lambda_1+\delta\lambda_2=0.\\
\end{array}\right.
\end{align}
Because $\alpha\neq0,$ ~~so we get $\lambda_2=\lambda_3.$~~Put~~$\lambda_2=\lambda_3$~~into ~~$\alpha^3\lambda_1+\frac{1}{2}(\lambda_3-\lambda_2)=0,$~~we have $\lambda_1=0.$~~So  another equation would reduce to~~$\delta\lambda_2=0,$~~and then $\delta=0.$~~\\
\end{pf}
\indent By (3.42) in \cite{C2}, we have
\vskip 0.5 true cm
\begin{lem}
Ricci symmetric tensors of $(G_7,g,J)$ ~~associated to Kobayashi-Nomizu connection are given by
\begin{align}
\widetilde{{\rm Ric}}^1\left(\begin{array}{c}
e_1\\
e_2\\
e_3
\end{array}\right)=\left(\begin{array}{ccc}
-\alpha^2&\frac{1}{2}(\beta\delta-\alpha\beta)&-\beta(\alpha+\delta)\\
\frac{1}{2}(\beta\delta-\alpha\beta)&-(\alpha^2+\beta^2+\beta\gamma)&-\frac{1}{2}(\beta\gamma+\alpha\delta+2\delta^2)\\
\beta(\alpha+\delta)&\frac{1}{2}(\beta\gamma+\alpha\delta+2\delta^2)&0
\end{array}\right)\left(\begin{array}{c}
e_1\\
e_2\\
e_3
\end{array}\right).
\end{align}
\end{lem}
By (4.23) and (2.3),we have
\begin{lem}
For $G_7$, the following equalities hold
\begin{align}
&{\rm L}_{e_1} \widetilde{\rm Ric}^1(e_1,e_2)=-\alpha^3+\frac{1}{2}\beta^2(3\delta+\alpha),{\rm L}_{e_1} \widetilde{\rm Ric}^1(e_1,e_3)=\alpha^3-\frac{1}{2}\beta^2(3\delta+\alpha),~~\\\notag
&{\rm L}_{e_1} \widetilde{\rm Ric}^1(e_2,e_2)=\beta(2\alpha\delta-3\alpha^2-2\beta^2-\beta\gamma+2\delta^2),{\rm L}_{e_1} \widetilde{\rm Ric}^1(e_2,e_3)=\beta(\frac{5}{2}\alpha^2+\frac{1}{2}\alpha\delta+\beta^2+\beta\gamma),~~\\\notag
&{\rm L}_{e_1} \widetilde{\rm Ric}^1(e_3,e_3)=-\beta(2\alpha^2+3\alpha\delta+\beta\gamma+2\delta^2),{\rm L}_{e_2} \widetilde{\rm Ric}^1(e_1,e_1)=2\alpha^2-\beta^2(3\delta+\alpha),~~\\\notag
&{\rm L}_{e_2} \widetilde{\rm Ric}^1(e_1,e_2)=\frac{1}{2}\beta(3\alpha^2+2\beta^2+\beta\gamma-2\alpha\delta-2\delta^2),{\rm L}_{e_2} \widetilde{\rm Ric}^1(e_1,e_3)=-\beta(\alpha^2+\frac{5}{2}\delta^2+2\alpha\delta+\frac{1}{2}\beta\gamma),~~\\\notag
&{\rm L}_{e_2} \widetilde{\rm Ric}^1(e_2,e_3)=\delta(\alpha^2+\beta^2-\delta^2-\frac{1}{2}\alpha\delta),{\rm L}_{e_2} \widetilde{\rm Ric}^1(e_3,e_3)=-\delta(3\beta\gamma+\alpha\delta+2\delta^2),~~\\\notag
&{\rm L}_{e_3} \widetilde{\rm Ric}^1(e_1,e_1)=-2\alpha^3+\beta^2(3\delta+\alpha),{\rm L}_{e_3} \widetilde{\rm Ric}^1(e_1,e_2)=\frac{1}{2}\beta(5\delta^2-3\alpha^2-2\beta^2-\beta\gamma+3\alpha\delta),~~\\\notag
&{\rm L}_{e_3} \widetilde{\rm Ric}^1(e_1,e_3)=\frac{1}{2}\beta(2\alpha^2+3\alpha\delta+\beta\gamma+2\delta^2),{\rm L}_{e_3} \widetilde{\rm Ric}^1(e_2,e_2)=\delta(2\delta^2+\alpha\delta-2\alpha^2-2\beta^2),~~\\\notag
&{\rm L}_{e_3} \widetilde{\rm Ric}^1(e_2,e_3)=\frac{1}{2}\delta(3\beta\gamma+\alpha\delta+2\delta^2).~~\\\notag
\end{align}
\end{lem}
By (4.24) and definition 2.3, we have
\vskip 0.2 true cm

If the unimodular Lorentzian Lie group $(G_7,g,J)$ is left-invariant Ricci collineations associated to the Kobayashi-Nomizu connection $\nabla^1$, then
${\rm L}_{\xi} \widetilde{\rm Ric}_{ij}^1=0,i=1,2,3~~ and ~~j=1,2,3$, so
\begin{align}
\left\{\begin{array}{l}
\left[2\alpha^2-\beta^2(3\delta+\alpha)\right]\lambda_2+\left[-2\alpha^3+\beta^2(3\delta+\alpha)\right]\lambda_3=0,\\
\\
\left[-\alpha^3+\frac{1}{2}\beta^2(3\delta+\alpha)\right]\lambda_1+\frac{1}{2}\beta(3\alpha^2+2\beta^2+\beta\gamma-2\alpha\delta-2\delta^2)\lambda_2+\frac{1}{2}\beta(5\delta^2-3\alpha^2-2\beta^2-\beta\gamma+3\alpha\delta)\lambda_3=0,\\
\\
\left[\alpha^3-\frac{1}{2}\beta^2(3\delta+\alpha)\right]\lambda_1-\beta(\alpha^2+\frac{5}{2}\delta^2+2\alpha\delta+\frac{1}{2}\beta\gamma)\lambda_2+\frac{1}{2}\beta(2\alpha^2+3\alpha\delta+\beta\gamma+2\delta^2)\lambda_3=0,\\
\\
\beta(2\alpha\delta-3\alpha^2-2\beta^2-\beta\gamma+2\delta^2)\lambda_1+\delta(2\delta^2+\alpha\delta-2\alpha^2-2\beta^2)\lambda_3=0,\\
\\
\beta(\frac{5}{2}\alpha^2+\frac{1}{2}\alpha\delta+\beta^2+\beta\gamma)\lambda_1+\delta(\alpha^2+\beta^2-\delta^2-\frac{1}{2}\alpha\gamma)\lambda_2+\frac{1}{2}\delta(3\beta\gamma+\alpha\delta+2\delta^2)\lambda_3=0,\\
\\
\beta(2\alpha^2+3\alpha\delta+\beta\gamma+2\delta^2)\lambda_1+\delta(3\beta\gamma+\alpha\delta+2\delta^2)\lambda_2=0.\\
\end{array}\right.
\end{align}
By solving (4.25), we get
\vskip 0.5 true cm
\begin{thm}
the unimodular Lorentzian Lie group $(G_7,g,J)$ is left-invariant Ricci collineations associated to the Kobayashi-Nomizu connection $\nabla^1$ if and only if\\
$(1)\alpha=0,\delta\neq0,\beta\neq0,\gamma=0,\lambda_2=\lambda_3=-\frac{\beta}{\delta}\lambda_1,\mathbb{V}_{\mathbb{R}\mathbb{C}}=\langle e_1-\frac{\beta}{\delta}e_2-\frac{\beta}{\delta} e_3\rangle,$\\
\\
$(2)\alpha=0,\delta\neq0,\beta=0,\lambda_2=\lambda_3=0,\mathbb{V}_{\mathbb{R}\mathbb{C}}=\langle e_1\rangle,$\\
\\
$(3)\alpha\neq0,\gamma=0,\delta=0,\lambda_1=0,\lambda_2=\lambda_3,\mathbb{V}_{\mathbb{R}\mathbb{C}}=\langle e_2+e_3\rangle.$\\
\end{thm}
\begin{pf}
  We analyze each one of these factors by separate,because $\alpha\gamma=0,$~~and if $\alpha=0$,~~we have~~$\delta\neq0$~~for $\alpha+\delta\neq0,$~~and Eq.(4.25) would reduce to
\begin{align}
\left\{\begin{array}{l}
\beta^2(\lambda_2-\lambda_3)=0,\\
\\
\frac{3}{2}\beta^2\delta\lambda_1+\frac{1}{2}\beta(2\beta^2\beta\gamma-2\delta^2)\lambda_2+\beta(\frac{5}{2}\delta^2-\beta^2-\frac{1}{2}\beta\gamma)\lambda_3=0,\\
\\
\frac{3}{2}\beta^2\delta\lambda_1+\beta(\frac{5}{2}\delta^2+\frac{1}{2}\beta\gamma)\lambda_2-\beta(\frac{1}{2}\beta\gamma+\delta^2)\lambda_3=0,\\
\\
\beta(-2\beta^2-\beta\gamma+2\delta^2)\lambda_1+\delta(2\delta^2-2\beta^2)\lambda_3=0,\\
\\
\beta(\beta^2+\beta\gamma)\lambda_1+\delta(\beta^2-\delta^2)\lambda_2+\delta(\frac{1}{2}\beta\gamma+\delta^2)\lambda_3=0,\\
\\
\beta(\beta\gamma+2\delta^2)\lambda_1+\delta(3\beta\gamma+2\delta^2)\lambda_2=0.\\
\end{array}\right.
\end{align}
If $\beta\neq0,$~~we have $\lambda_2=\lambda_3.$~~Put $\lambda_2=\lambda_3$~~into Eq.(4.26), we get
\begin{align}
\left\{\begin{array}{l}
\beta\lambda_1+\delta\lambda_3=0,\\
\\
\gamma\lambda_1=0.\\
\end{array}\right.
\end{align}
So $\gamma=0,-\frac{\beta}{\delta}\lambda_1=\lambda_2=\lambda_3$,~~otherwise $\gamma\neq0,\lambda_1=\lambda_2=\lambda_3=0.$\\
If $\beta=0,$~~then we get
\begin{align}
\left\{\begin{array}{l}
\delta^3\lambda_3=0,\\
\\
\delta^3\lambda_2=0,\\
\\
\delta^3(\lambda_3-\lambda_2)=0.\\
\end{array}\right.
\end{align}
So we have $\lambda_2=\lambda_3=0,\mathbb{V}_{\mathbb{R}\mathbb{C}}=\langle e_1\rangle.$\\
If $\alpha\neq0,\gamma=0,$~~ we get
\begin{align}
\left\{\begin{array}{l}
\left[2\alpha^2-\beta^2(3\delta+\alpha)\right](\lambda_2-\lambda_3)=0,\\
\\
\left[-2\alpha^2+\beta^2(3\delta+\alpha)\right]\lambda_1+\beta(3\alpha^2+2\beta^2-2\alpha\delta-\delta^2)\lambda_2+\beta(5\delta^2-3\delta^2-2\beta^2+3\alpha\delta)\lambda_3=0,\\
\\
\left[-2\alpha^2+\beta^2(3\delta+\alpha)\right]\lambda_1+\beta(2\alpha^2+4\alpha\delta+5\delta^2)\lambda_2-\beta(2\alpha^2+3\alpha\delta+2\delta^2)\lambda_3=0,\\
\\
\beta(\frac{5}{2}\alpha^2+\frac{1}{2}\alpha\delta+\beta^2)\lambda_1+\delta(\alpha^2+\beta^2-\delta^2-\frac{1}{2}\alpha\delta)\lambda_2+\delta(\frac{1}{2}\alpha\delta+\delta^2)\lambda_3=0,\\
\\
\beta(2\alpha^2+3\alpha\delta+2\delta^2)\lambda_1+\delta(\alpha\delta+2\delta^2)\lambda_2=0.\\
\end{array}\right.
\end{align}
If $2\alpha^2-\beta^2(3\delta+\alpha)\neq0,$~~then we have $\lambda_2=\lambda_3.$~~Put $\lambda_2=\lambda_3$~~into other equations, we get
\begin{align}
\left\{\begin{array}{l}
\left[-2\alpha^2+\beta^2(3\delta+\alpha)\right]\lambda_1+\beta(3\delta^2+3\alpha\delta)\lambda_3=0,\\
\\
\beta(\frac{5}{2}\alpha^2+\frac{1}{2}\alpha\delta+\beta^2)\lambda_1+\delta(\alpha^2+\beta^2)\lambda_3=0,\\
\\
\beta(2\alpha^2+3\alpha\delta+2\delta^2)\lambda_1+\delta(\alpha\delta+2\delta^2)\lambda_3=0.\\
\end{array}\right.
\end{align}
When $\delta=0,$~~so $2\alpha^2-\beta^2\neq0.$~~Put this into other equations, Eq.(4.30)would reduce to~~$\lambda_1=0$, so $\mathbb{V}_{\mathbb{R}\mathbb{C}}=\langle e_2+e_3\rangle.$~~\\
When $\delta\neq0,$~~we get $\lambda_1=\lambda_2=\lambda_3=0.$~~ \\
If  $2\alpha^2-\beta^2(3\delta+\alpha)=0,$~~$\beta\neq0,$~~otherwise $\alpha=0$~~in error.\\Furthermore if $\delta=0,$~~we get $2\alpha^2-\beta^2=0,\lambda_0.$~~Put these into Eq.(4.30)~~, we get $\lambda_2=\lambda_3.$\\
If $\delta\neq0,$~~then we have $\delta=\frac{\alpha(2\alpha^2-\beta^2)}{3\beta^2}.$~~ We assume that
\begin{align}
\left\{\begin{array}{l}
a_{11}=0,a_{12}=-(\alpha^2+\frac{5}{2}\delta^2+2\alpha\delta),a_{13}=\alpha^2+\frac{3}{2}\alpha\delta+\delta^2\\
\\
a_{21}=\frac{5}{2}\alpha^2+\frac{1}{2}\alpha\delta+\beta^2,a_{22}=\delta(\alpha^2+\beta^2-\delta^2-\frac{1}{2}\alpha\delta),a_{23}=\delta(\frac{1}{2}\alpha\delta+\delta^2)\\
\\
a{31}=\beta(2\alpha^2+3\alpha\delta+2\delta^2),a_{32}=\delta(\alpha\delta+2\delta^2),a_{33}=0.\\
\end{array}\right.
\end{align}
We get
\begin{align}
\left|\begin{array}{ccc}
a_{11}&a_{12}&a_{13}\\
a_{21}&a_{22}&a_{23}\\
a_{31}&a_{32}&a_{33}\\
\end{array}\right|
\neq0,
\end{align}
So we have~~$\lambda_1=\lambda_2=\lambda_3=0.$~~
\end{pf}
\vskip 1 true cm

\section{Acknowledgements}

The author are deeply grateful to the referees for their valuable commments and helpful suggestions.

\vskip 1 true cm


\bigskip
\bigskip

\noindent {\footnotesize {\it Tao Yu} \\
{School of Mathematics and Statistics, Northeast Normal University, Changchun 130024, China}\\
{Email: yut338@nenu.edu.cn}\\

\end{document}